\documentclass[12pt,a4paper]{article}
\usepackage[utf8]{inputenc}
\usepackage[T1]{fontenc}
\usepackage{lmodern}

\usepackage{float}
\usepackage{amsmath,amssymb,amsthm,mathtools}
\usepackage{bm}
\usepackage{multirow}
\usepackage{enumitem}
\usepackage{graphicx}
\graphicspath{{figures/}}
\usepackage{caption}
\usepackage{subcaption}
\usepackage[colorlinks=true,linkcolor=blue,citecolor=blue,urlcolor=blue]{hyperref}
\usepackage{algorithm}
\usepackage{algpseudocode}
\usepackage{booktabs}
\usepackage{geometry}
\newtheorem{proposition}{Proposition}[section]
\newtheorem{theorem}{Theorem}[section]
\newtheorem{lemma}{Lemma}[section]
\theoremstyle{definition}
\newtheorem{definition}{Definition}[section]
\newtheorem{remark}{Remark}[section]
\newtheorem{example}{Example}[section]
\newtheorem{corollary}{Corollary}[section]

\newcommand{\precprec}{\prec\mathrel{\mkern-5mu}\prec}
\newcommand{\IR}{\mathbb{IR}}
\newcommand{\R}{\mathbb{R}}

\newcommand{\prox}{\operatorname{prox}}

\title{\textbf{Proximal point algorithm for interval valued optimization problems with application to lasso problem }}
\author{
Anurag Jayswal\textsuperscript{a},
   Sunidhi\textsuperscript{{*}a}, 
   Caihong Han\textsuperscript{b}
   }
\date{}

\begin{document}
\maketitle
\renewcommand{\thefootnote}{*}
\footnotetext{%
Corresponding author\\
E-mail addresses: anurag@iitism.ac.in (A. Jayswal), *sunidhis2031@gmail.com (Sunidhi),\\
caihonghansju@163.com (Caihong Han)
}
\begin{center}
   \textit{\textsuperscript{a}Department of Mathematics and Computing, Indian Institute of Technology (Indian School of Mines), Dhanbad-826004, Jharkhand, India}\\
   \textit{\textsuperscript{b}School of Mathematics and Quantitative Economy, Guangxi University of Finance and Economics, Nanning, 53003, China}
\end{center}

\section*{Abstract}
\noindent In this paper, we consider a class of unconstrained interval valued optimization problems and using a suitable scalarization technique, we formulate an associated scalarized optimization problem and establish the relationships between the optimal solution of the scalarized problem and the $LU$-nondominated solution of the considered  optimization problem. Based on this framework, we propose a proximal point algorithm and derive the global convergence of the sequence generated by the proposed algorithm along with the sublinear and linear rate of convergence for convex and strongly convex cases, respectively. Finally, we provide suitable numerical experiments along with comparative analysis with an existing method and applicability of the proposed algorithm on an interval valued $lasso$ regression problem, highlighting its potential in machine learning models with uncertain data.

\medskip
\noindent \textbf{Keywords:}
Interval valued optimization $\cdot$ Scalarization $\cdot$   Proximal point algorithm $\cdot$   $LU$-ordering.

\medskip
\noindent \textbf{Mathematics Subject Classification} 90C25 $\cdot$ 90C70 $\cdot$ 90C29 $\cdot$ 65K10 $\cdot$ 65G40

\section{Introduction}

Optimization problems with uncertain or imprecisely known data arise naturally in engineering,
economics, machine learning, and decision-making, etc. Interval analysis, introduced by Moore
\cite{moore1966interval} and further consolidated by Moore $et\ al.$\ \cite{Moore2009} and
Alefeld and Herzberger \cite{Alefeld1983}, provides a rigorous mathematical framework for handling such
uncertainty by representing each uncertain quantity as a closed bounded interval. An
interval-valued optimization problem (IVOP) arises when the objective function
maps decision variables to intervals rather than real numbers. A central difficulty in (IVOP) is the absence of a total ordering on the interval space,
which necessitates the use of partial orderings to define appropriate solution concepts.
The $LU$-ordering, introduced by Ishibuchi and Tanaka \cite{Ishibuchi1990}, compares
intervals componentwise via their lower and upper endpoints and has become the standard
framework for defining LU-nondominated solutions, in analogy with Pareto
optimality in multiobjective optimization \cite{Ehrgott2005}. Building on this, Wu
\cite{Wu2007,Wu2008} established Karush--Kuhn--Tucker (KKT) optimality conditions for (IVOP)
under LU-ordering, while duality and saddle point results were further studied by
Jayswal $et\ al.$\ \cite{Jayswal2011} and Ahmad $et\ al.$\ \cite{Ahmad2016}.\\

Parallel to the optimality theory, a variety of solution methods have been proposed for (IVOP), most of them built on the notion of the generalized Hukuhara ($gH$) difference introduced by Stefanini~\cite{Stefanini2009} to overcome the limitations of the classical Hukuhara difference \cite{Hukuhara1967}. The $gH$-difference underlies the $gH$-differentiability framework developed by Stefanini and Bede \cite{Stefanini2010} and Chalco-Cano $et\ al.$\cite{Chalco2011}, and has enabled the design of a range of gradient-based and Newton-type algorithms for smooth (IVOP). Among the Newton-type methods, Ghosh \cite{ghosh2017newton} proposed a Newton method to obtain the efficient solutions of (IVOP) with $gH$-differentiable interval valued objectives which was later extended to multiobjective case by Upadhyay $et\ al.$~\cite{upadhyay2024newton}. Further in the gradient-based approaches, Ghosh $et\ al.$~\cite{ghosh2022generalized} introduced a generalized-Hukuhara gradient efficient direction method for solving interval valued least squares problems and derived the optimality conditions of the (IVOP) using a novel definition of $gH$-differentiability. Recently, Mondal and Ghosh~\cite{mondal2025steepest} developed a steepest descent method for multiobjective interval-valued optimization problem and characterized the connections between the
weakly Pareto optimal points and Pareto stationary points.\\

To address the class of nonsmooth (IVOP), Liu $et\ al.$~\cite{Liu2023} developed a subgradient-based neurodynamic algorithm for constrained nonsmooth, nonconvex (IVOP) involving partial-order together with linear equality constraints, formulating the iterative scheme as a differential inclusion. Further, Guo $et\ al.$~\cite{Guo2023} employed invexity-type generalized convexity notions to handle nonsmooth (IVOP) and derived the optimality and duality results using the symmetric invexity framework. Recently, Kumar $et.al$~\cite{KSR} proposed a projected $gH$-subgradient method for the class of nonsmooth (IVOP). While these approaches have proved effective in computing various classes of (IVOP), they generally rely on some form of differentiability or at least $gH$-subdifferentiability of the interval valued objective function, thereby limiting their applicability.\\

Beyond gradient-based schemes, the proximal point algorithm offers a fundamentally different route for handling nonsmooth convex minimization problems, since it relies only on convexity and lower semicontinuity of the objective rather than any differentiability assumption. Originally introduced by Martinet~\cite{Martinet1970} as a regularization scheme for variational inequalities, it was extensively developed by Rockafellar \cite{rockafellar1976monotone, rockafellar1970maximal} for maximally monotone operators and subdifferential mappings. The convergence behavior of the proximal point algorithm for the nonsmooth convex minimization problems has been established through the classical variational analysis of Rockafellar and Wets~\cite{rockafellar1998variational} and Rockafellar~\cite{tyrrell1970convex} which was further refined by Güler \cite{Guler1991}, Kamimura and Takahashi \cite{Kamimura2003}, Parikh and Boyd \cite{parikh2014proximal} and Beck \cite{Beck2017}.The proximal framework also underpins a broad family of modern first-order methods, including proximal gradient method, operator splitting techniques, and the alternating direction method of multipliers \cite{beck2009fast, beck2010gradient, combettes2005signal,Boyd2011, Combettes2011}.
Among the literature survey consolidating these developments, we particularly draw on and refer to the comprehensive analysis of the proximal point algorithm by Cai $et\  al.$~\cite{cai2022developments}, which unifies the classical, inexact and accelerated variants of the proximal point algorithm along with their convergence and rate analysis. Although the proximal point theory has already been extensively investigated in the classical real-valued setting, to the best of our knowledge no work to date has introduced a proximal point algorithm to address the class of interval valued optimization problems.\\

Motivated by this gap, we propose a proximal point algorithm for solving a class of interval valued optimization problem (IVOP) whose objective functions are proper, convex, and generalized Hukuhara lower semicontinuous. The (IVOP) is subsequently transformed into a scalarized optimization problem (SOP) through a suitable scalarization technique. We then establish that properness, convexity and $gH$-lower semicontinuity of the interval-valued objective are exactly the conditions needed to guarantee that the associated scalarized objective is well posed, thereby ensuring the existence of a unique solution to every proximal subproblem. The main contributions of this paper are as follows:
\begin{itemize}
\item We establish a relationship between the $LU$-nondominated solution of the (IVOP) and the optimal solution of the associated (SOP).
  \item We prove that the associated scalarized objective inherits properness, convexity, lower semicontinuity and strong convexity from the interval valued objective. We further establish a sufficient condition under which an $LU$-nondominated solution of the (IVOP) can be recovered from the scalarized objective through $gH$-subdifferential characterization.

  \item Based on this framework, we propose the proximal point algorithm for the interval valued optimization problem, referred as PPA-IVOP, and prove that each proximal subproblem is well defined and admits a unique minimizer at every iteration.

  \item We establish the descent property of the sequence generated by the PPA-IVOP algorithm and prove that every cluster point of the sequence of iterates is an $LU$-nondominated solution of the (IVOP). We further provide sublinear convergence under convexity and linear convergence under strong convexity.
  \item We identify a structural condition namely, the existence of a common minimizer of the lower and upper endpoint functions of the objective, under which the proposed algorithm converges to a unique $LU$-nondominated solution, independent of the choice of the scalarization parameter.
  \item Numerical experiments along with a comparative analysis are presented to validate the theoretical results and effectiveness of the proposed algorithm and finally, we demonstrate the applicability of the algorithm on an interval-valued $lasso$ regression problem.
\end{itemize}

The remainder of this paper is organized as follows. Section~\ref{Sect.2} reviews the necessary preliminaries from proximal and interval analysis. Section~\ref{Section 3} formulates the (IVOP) together with its associated (SOP) and establishes the scalarization framework underlying the proposed algorithm and presents the PPA-IVOP algorithm. Section~\ref{Section Convergence} presents the convergence analysis of the proposed algorithm
and the numerical experiments to illustrate its performance together with a comparative analysis. Section~\ref{Section Application} demonstrates the applicability of the proposed algorithm to an interval valued $lasso$ regression problem. The paper concludes in Section~\ref{Section Conclusion} by summarizing the key results and highlighting promising directions for future studies.

\section{Preliminaries}\label{Sect.2}

In this section, we begin with some standard notations used throughout the article. Next, we recall some basic definitions and concepts from proximal analysis that are fundamental tool of the proposed methodology. Finally, we present the necessary notions from the interval analysis that will be used in the subsequent sections.\\
\\
\noindent We adopt the following notations in this article:  
$\mathbb{R}$, $\R_+$, $\overline{\mathbb{R}} = \mathbb{R} \cup \{ +\infty\}$,
$\mathbb{R}^n$,
$\IR$,
$\overline{\IR} = \IR \cup \{[ +\infty,+\infty]\}$, $\IR_+$ and $\IR^n$ to represent the real line, the non-negative real line, the extended real line, $n$-dimensional Euclidean space, the set of compact intervals defined on the real line, the family of compact intervals defined on the extended real line, the family of compact intervals on non-negative real line and the Cartesian product space of compact intervals $\IR\times \IR\times \dots\times\IR(n$ times$)$, respectively.

\subsection*{Proximal Analysis}\label{prox anlys}

Let $\phi:\R^n\to \overline{\R}$ be an extended real valued function. 

\begin{definition}[\textit{Closed function}] The function $\phi$ is termed as closed if its epigraph, given by
\[
epi \ \phi=\big\{(x,t)\in \R^{n+1} : \phi(x)\le t\big\},
\]
 is closed which means every sequence $\{(x_k,t_k)\}\subset epi\phi$ converges to some $(x,t)\in epi\phi$.   
\end{definition}

\begin{definition}[\textit{Effective domain}]\label{eff dom scalar} The effective domain of the function $\phi$, denoted by $dom\ \phi$, is given as
\[
dom\ \phi=\big\{ x\in \R^n : \phi(x)<+\infty \big\}.
\]
\end{definition}

\begin{definition}[\textit{Proper function} \cite{tyrrell1970convex}]\label{def:proper func} The function $\phi$ is said to be proper if:
\begin{enumerate}[label=(\roman*)]
    \item $\phi(x)>-\infty, \ \forall x\in \R^n$,
    \item there exists at least one $x_0\in \R^n$ such that $\phi(x_0)<+\infty$.
\end{enumerate}
\end{definition}

\begin{definition}[\textit{Convex function} \cite{tyrrell1970convex}] The function $\phi$ is said to be convex if for all $x,y \in dom\ \phi$ and $\beta \in [0,1]$, one has
\[
\phi(\beta x + (1-\beta)y)
\le
\beta \phi(x) + (1-\beta)\phi(y).
\]
If the above inequality holds with equality for $\beta\in[0,1]$, then $\phi$ is called an affine function, for all $x,y\in dom \ \phi$.
    
\end{definition}

\begin{definition}[\textit{Strong Convexity} \cite{tyrrell1970convex}]\label{def:strongly convex}
 The function $\phi$ is said to be
$\gamma$-strongly convex, with modulus $\gamma > 0$, if
\[
\phi(\beta x + (1-\beta)y)
\le
\beta \phi(x) + (1-\beta)\phi(y)
- \frac{\gamma}{2}\beta(1-\beta)\|x-y\|^2,
\]
for all $x,y \in dom\ \phi$ and $\beta \in [0,1]$.
\end{definition}

\begin{definition}[\textit{Lower Semicontinuity}]\label{lsc} The function $\phi$ is said to be lower semi-continuous ($l.sc$) at a point $x_0\in \R^n$, if we have
\[
   \liminf_{x \to x_0} \phi(x) \ge \phi(x_0), \]
where
\[
\liminf_{x \to x_0} \phi(x)=\sup_{\delta > 0}\inf\{\phi(x):x\in \R^n, 0<\|x-x_0\|<\delta\}.
\]
\end{definition}

\begin{definition}[\textit{Normal Cone} \cite{tyrrell1970convex}]\label{ncone} Let $\mathcal T$ be a nonempty convex subset of $\R^n$. The normal cone to $\mathcal T$ at $x_0 \in \mathcal T$, denoted by $N_{\mathcal T}(x_0)$, is defined as
\[
N_{\mathcal T}(x_0)
:=
\left\{
y \in \mathbb{R}^n \;\middle|\;
\langle  y, x-x_0 \rangle \le 0,
\ \forall x \in \mathcal T
\right\}.
\]
\end{definition}

\begin{definition}[\textit{Coercive Function}]
The function $\phi $
is said to be coercive if for every sequence $\{x^k\}\subset \R^n$ satisfying $\|x^k\|\to+\infty$, we have $\phi(x_k)\to +\infty$.
\end{definition}

\begin{definition}[\textit{Subdifferential}\label{def: subdiff}] 
 Let $\phi : \R^n\to \overline{\mathbb{R}}$ be a convex function. The subdifferential of $\phi$ at a point
$x^* \in dom\ \phi$ is defined as
\begin{align*}\label{sugrad def eq}
\partial \phi(x^*)
:=
\left\{
g \in \mathbb{R}^n \;\middle|\;
\phi(x) \ge \phi(x^*) + \langle g, x - x^* \rangle,
\ \forall x \in \R^n
\right\},
\end{align*}
where $g$ is called the subgradient of the function $\phi$ at the point $x^*$. If $x^* \notin dom\ \phi$, then $\partial \phi(x^*) = \emptyset$.

\end{definition}

\begin{remark}\label{sugdiff remark}
If the function $\phi$ is proper convex and $x\in int(dom\ \phi)$, then its subdifferential $\partial \phi(x)$ is a nonempty closed convex set. Further, if the function $\phi$ is $\gamma$-strongly convex, then $\phi$ satisfies the following inequality
\[
\phi(x)\ge \phi(x^*) +\langle g, x - x^* \rangle+\frac{\gamma}{2}\|x-x^*\|^2, 
\]
for all $x\in \R^n$, where $g\in \partial\phi(x^*)$ and $x^*\in dom \ \phi$.
\end{remark}

\begin{definition}[\textit{Proximal operator}\cite{parikh2014proximal}] Let $\phi:\R^n\to\overline{\R}$ be a proper convex, closed and lower semicontinuous function, then the proximal operator $\prox_{\mu_k\phi}:\R^{n}\rightarrow\mathcal \R^n$ of $\phi$ with parameter $\mu_k>0$ is defined by
\[
\operatorname{prox}_{\mu_k \phi}(y)
=
\arg\min_{x\in\R^n}
\left\{
\phi(x)+\frac{1}{2\mu_k}\|x-y\|^2
\right\}, \ \forall y\in\R^{n}.
\]
Here, the lower semicontinuity of the function $\phi$ ensures the existence of minimizers of the proximal subproblem and guarantees that the proximal operator is well defined.
\end{definition}
\begin{remark}
    Since the function $\phi$ is proper, convex and lower semicontinuous, the function $x\mapsto \phi(x)+\frac{1}{2\mu_k}\|x-y\|^2\}$ is proper and strongly convex for every $y\in\R^n$, and hence it admits a unique minimizer for every $y\in \R^n$. 
Equivalently, one can say that the proximal operator finds a point $x\in \R^n$ that minimizes the value of $\phi(x)$ while staying close to the point $y$.
\end{remark}

\subsection*{Interval Analysis}
Let $G\in \IR$, then $G$ can be represented as $G=[G^l,G^u]$, where $G^l\le G^u$ and $G^l$,$G^u\in \R$.
According to Moore \cite{moore1966interval}, for a pair of intervals
$G=[G^l,G^u]$ and $H=[H^l,H^u]$ in $\IR$, and a real number $c$, the interval arithmetic are denoted as $G\oplus H$, $G\ominus H$, $c\odot G$, and $G\odot H$, respectively, and are defined as:
\begin{itemize}
    \item $G\oplus H=[G^l+H^l,\ G^u+H^u].$
    \item $G\ominus H=[G^l-H^l,\ G^u-H^u].$
    \item $ c \odot G =
    \begin{cases}
[c G^l,\ c G^u], & c \ge 0, \\
[c G^u,\ c G^l], & c < 0 .
\end{cases}$
   \item $ G\odot H$= [min$\{G^lH^l, G^lH^u, G^uH^l,G^uH^u\}$, max$\{G^lH^l, G^lH^u, G^uH^l,G^uH^u\}].$
\end{itemize}

In addition, if $G\in \IR$ is such that $G^l=G^u$, then $G$ is called a degenerate interval. In this case, we denote
\[
\textbf{G}=[v,v],\quad v\in \R,
\]
where $v=G^l=G^u$. In particular, the zero interval is denoted by $\textbf{0}=[0,0]$.\\
In interval analysis, order relations are used to compare the interval values and thus, several partial orderings of the intervals have been discussed in the literature (see \cite{moore1966interval, ghosh2017newton, wu2011interval, mondal2025steepest}, and the references cited therein). We consider the partial ordering between the intervals $G,H \in \IR$ from Wu \cite{wu2011interval} as :
\[
G \preceq_{LU} H \quad \Leftrightarrow \quad G^l\le H^l
\ \text{and}\ 
G^u\le H^u,
\]
where "$\preceq_{LU}$" refers to the partial ordering in $\IR$. Thus, $G$ is inferior to $H$ if $G$ lies entirely inside $H$. Also, $G\prec_{LU}H$ if and only if $G\preceq_{LU}H$ and $G\neq H$ which means $G\prec_{LU}H$ if and only if :
\[
\begin{cases}
    G^l< H^l\\
    G^u=H^u
\end{cases}
,\quad \text{or} \quad
\begin{cases}
    G^l= H^l\\
    G^u<H^u
\end{cases}
,\quad \text{or} \quad\
\begin{cases}
    G^l< H^l\\
    G^u<H^u
\end{cases}
.
\]
We also write $G\precprec_{LU} H$ if and only if $G^l<H^l$ and $G^u<H^u$.\\
In the interval-valued arithmetic, the usual difference of two compact intervals is not suitable for interval calculus, since it does not preserve basic algebraic properties (see Moore \cite{moore1966interval}). To overcome this challenge, Stefanini \cite{stefanini2008generalization} introduced the concept of generalized Hukuhara difference.

\begin{definition}[\textit{Generalized Hukuhara ($gH$)-difference}\cite{stefanini2008generalization}]\label{def: gh difference} Let $G,H$ and $J$ be three arbitary elements of $\IR$. If $G=H\oplus J$ or $H=G\ominus J$, then $J$ is called the \textit{$gH$-difference} between $G$ and $H$, and is denoted as $J=G\ominus_{gH} H$. For a given pair of intervals $G=[G^l,G^u]$ and $H=[H^l,H^u]$, the $gH$-difference is defined as:
\[
   G\ominus_{gH}H=[min\{G^l-H^l,\ G^u-H^u\},\ max\{G^l-H^l,\ G^u-H^u\}].
\]
\end{definition}

\begin{definition}[\textit{Norm on $\IR$}]\label{IVF norm}
Let $G=[G^l,G^u]\in \IR$, then the norm on $\IR$ is given by
\[
   \|G\|_{\IR}=max\{|G^l|,|G^u|\}.
\] 
\end{definition}

\subsection*{Interval valued functions}
\begin{definition}[\textit{Interval valued function}]\label{def:IVF} A function $F : \R^n \to \IR$ is called an interval valued function (IVF), if for every $x\in \R^n$,
\[
F(x) = [F^l(x), F^u(x)],
\]
where the endpoint functions $F^l, F^u : \R^n \to \R$ satisfy
$F^l(x) \le F^u(x), \ \forall x \in \R^n$.
\end{definition}

\begin{remark}\label{remark: ext ivf}
    If $F$ takes values in the extended interval space $\overline{\IR}$, then $F$ is called an extended interval-valued function. In this setting, the endpoint functions of $F$ take values in the extended real line $\overline{\R}$
\end{remark}

\begin{definition}[\textit{$gH$-derivative}] Let $F:\mathcal T\to \IR$ be an interval valued function defined on a non-empty open convex set $\mathcal T\subseteq\R^n$. The generalized Hukuhara derivative of $F$ at a point $x^*\in \mathcal T$, denoted by $D_{gH}F(x^*)$, is defined as 
\[
   D_{gH}F(x^*)=\lim_{d\to 0} \frac{1}{d}\odot \{F(x^*+d)\ominus_{gH}F(x^*)\}.
\]
 If the above limit exists, then the IVF $F$ is said to be generalized Hukuhara differentiable ($gH$-differentiable) at $x^*\in \mathcal T$. 
\end{definition}

\begin{definition}[\textit{Effective domain of IVF} \cite{anshika2022interval}]\label{effective domain ivf def} The effective domain of an extended IVF $F: \R^n\rightarrow \overline{\IR} $ is defined as:
\[
  dom\ F=\{x\in \R^n:\|F(x)\|_{\IR}<\infty\}.
\]
\end{definition}

\begin{definition}[\textit{Proper IVF} \cite{kumar2023ekeland}]\label{def:proper}
An extended IVF $F: \R^n\rightarrow \overline{\IR} $ is said to be proper if the following properties hold:
\begin{enumerate}[label=(\roman*)]
 \item $[-\infty,-\infty]\precprec_{LU} F(x), \ \forall x\in \R^n$,
 \item $F(x)\precprec_{LU}[+\infty,+\infty]$ for at least one $x\in \R^n$,
\end{enumerate}
\end{definition}

\begin{remark} \label{Remark 2.4}
In the view of Definitions \ref{def:proper} and \ref{def:proper func}, it follows that if the IVF $F$ is proper, then the endpoint functions are also proper.
 \end{remark}

\begin{definition}[\textit{Convex IVF}\label{def:convex}] An extended IVF $F:\R^n \rightarrow \overline{\IR}$, is called convex, if 
\[
F(\beta x + (1-\beta)y)
\preceq_{LU}
\beta\odot F(x) \oplus (1-\beta)\odot F(y),
\]
for all $x,y\in dom\ F$
and $\beta\in[0,1]$.
Equivalently,
\[
F^l(\beta x + (1-\beta)y)
\le
\beta F^l(x) + (1-\beta)F^l(y),
\]
\[
F^u(\beta x + (1-\beta)y)
\le
\beta F^u(x) + (1-\beta)F^u(y),
\]
\end{definition}

\begin{definition}[\textit{Strongly convex IVF}\label{def:strgly convex}] An extended IVF $F:\R^n \rightarrow \overline{\IR}$, is called $\gamma$-\textit{strongly convex} whenever there exists a convex interval valued function $H:\R^n\to \overline{\IR}$ and a constant $\gamma>0$ satisfying
\[
F(x)=H(x)\oplus{\gamma}\odot\mathbf{\frac{\|x\|^2}{2}},\quad  \forall x\in dom\ F.
\]
\end{definition}

\begin{remark}\label{remark convex}
   Definitions \ref{def:convex} and \ref{def:strgly convex} immediatly imply that the IVF $F$ is (strongly) convex if and only if $F^l$ and $F^u$ are (strongly) convex functions (see for instance \cite{ghosh2022generalized}).
\end{remark}

\begin{definition}[\textit{Infimum and supremum of an extended IVF} \cite{kumar2023ekeland}]For a nonempty set $\mathcal T\subseteq\R^n$, consider an extended IVF $F:\mathcal T\to \overline{\IR}$. The \emph{infimum} of $F$ over $\mathcal T$ is defined as the infimum of its range set, namely
\[
\inf_{x \in \mathcal T} F(x)
:=
\inf\{F(x):x\in\mathcal{T}\}.
\]
In the analogous manner, the \emph{supremum} of $F$ over $\mathcal T$ is given by
\[
\sup_{x \in \mathcal T} F(x)
:=
\sup\{F(x):x\in\mathcal{T}\}.
\]
\end{definition}

\begin{definition}[\textit{$gH$-Lower Semicontinuity of an extended IVF} \cite{kumar2023ekeland}]\label{def:gHLsc}Let $F:\R^n \to \overline{\IR}$ be an extended IVF and let $x^*\in \R^n$. For $\delta>0$, we define the open ball $B_{\delta}(x^*):=\{x\in \R^n: \|x-x^*\|<\delta\},$ and the limit inferior of $F$ at $x^*$ as
\[
\lim_{x\to x ^*}\inf F(x)=\sup_{\delta >0}\inf \{F(x): x\in B_\delta (x^*)\}.
\]
Then, the IVF $F$ is a generalized Hukuhara lower semicontinuous ($gH$-$l.sc$) at $x^*\in \R^n$, if
\[
F(x^*)\preceq_{LU}\lim_{x\to x ^*}\inf F(x),\quad \forall x\in B_{\delta}(x^*).
\]
Further, the IVF $F$ is said to be $gH$-$l.sc$ on $\R^n$, if the above inequality holds for every $x^*\in \R^n$.
\end{definition}
\begin{remark}\label{Remark 2.5} In view of the Definition \ref{def:gHLsc}, it follows that the extended IVF $F$ is $gH$-lower semicontinuous at the point $x^*\in \R^n$ if and only if the end point functions are lower semicontinuous at the point $x^*\in \R^n$ in the sense of Definition \ref{lsc}. (Refer to Note $3.1$ of Kumar and Ghosh \cite{kumar2023ekeland})
\end{remark}

\begin{definition}[\textit{$gH$-continuity}]\label{gh continuity} 
 An extended IVF $F:\R^n\to \overline{\IR}$ is said to be $gH$-continuous a point $x^*\in\R^n $, provided that
\[
\lim_{x\to x^*}\|F(x)\ominus_{gH}F(x^*)\|_{\IR}=0.
\]
If the above condition holds for every $x^*\in \R^n$, then the IVF $F$ is $gH$-continuous on $\R^n$.
    
\end{definition}

\begin{remark}\label{remark continuity}
    In the view of Definition~\ref{gh continuity}, it follows that the extended IVF $F$ is $gH$-continuous on $\R^n$ if and only if the endpoint functions are continuous on $\R^n$. Furthermore, if an extended IVF $F$ is $gH$-continuous on $\R^n$, then $F$ is $gH$-lower semicontinuous on $\R^n$ (see Kumar and Ghosh \cite{kumar2023ekeland}).
\end{remark}

\begin{definition}[\textit{$gH$-subdifferential}]\label{gh-subdiff}
    Let $F:\R^n\to \overline{\IR}$ be a proper convex IVF. Then the $gH$-subdifferential of $F$ at a point $x_0\in dom\ F$, is the set-valued map given by
    \[
    \partial _{gH}F(x_0)=\{\hat{\mathcal G}\in \IR^n | \ (x-x_0)^\top\odot [g^-,g^+]\preceq_{LU} F(x)\ominus_{gH}F(x_0),\ \forall x\in \R^n\},
    \]
    where $\hat{\mathcal G}=[g^-,g^+]$ with $g^- =(g^- _1,g^- _2,\dots,g^- _n),\ g^+ =(g^+ _1,g^+ _2,\dots,g^+ _n)\in \R^n$. And, if $\partial _{gH}F(x_0)\neq \emptyset$, then $F$ is said to be $gH$-subdifferentiable at $x_0\in dom\ F$.
\end{definition}

\begin{remark}\label{gh-subdiff remark}
By componentwise analysis, 
$\mathcal {\hat G}\in\partial_{gH}F(x_0)$ if and only if
$g^-\in\partial{F}^l(x_0)$ and
$g^+\in\partial {F}^u(x_0)$, (see, for instance, Kumar $et\ al.$ \cite{Kumar2022weak}).
\end{remark}

\begin{theorem}\normalfont{\cite{Kumar2022weak}}
\label{lem:subdiff_nonempty}
If the IVF $F:\R^n\to \overline{\IR}$ is a proper convex function and $x_0\in int(dom\ F)$, then
$\partial_{gH}F(x_0)$ is a closed and convex set.
\end{theorem}

\section{Proximal Point Algorithm for Interval-Valued\\ Optimization}\label{Section 3}

In this section, we consider an unconstrained interval-valued optimization problem (IVOP) and formulate an associated scalarized optimization problem (SOP). Further, the relationship between the solution of the two problems is established along with the preservation of the structural properties of the interval valued objective. Based on this framework, a proximal point algorithm is proposed to compute solutions of the considered (IVOP).
\\
\\
Consider the unconstrained interval-valued optimization problem, defined as follows:
\begin{equation*}
    \textbf{(IVOP)} \qquad \min_{x\in \R^n} F(x), \label{IVOP}
\end{equation*}
where $F:\R^n \to \overline{\IR}$ is a proper convex and $gH$-lower semicontinuous extended IVF. Under these assumptions, we call the
IVF $F$ a PCL function (Proper, Convex and $gH$-Lower semicontinuous IVF). 
\begin{definition}
    Let $x^*\in \R^n$. Then $x^*$ is referred to an \emph{$LU$-nondominated solution}  of (IVOP), if there does not exist any other point $x\in \R^n$, such that
    \begin{equation*}
        F(x)\prec_{LU}F(x^*).
    \end{equation*}
\end{definition}

Now, we formulate a scalarized optimization problem (SOP), associated with the (IVOP), in the following manner:
\begin{equation*}
    \textbf{(SOP)} \qquad \min_{x\in \R^n} \phi_\lambda(x), \label{SOP}
\end{equation*}
where $\phi_\lambda:\R^n \to \overline{\R}$, is defined as:
\begin{equation}\label{scalar function}
\phi_\lambda(x):=\lambda F^l(x)+(1-\lambda)F^u(x), \quad \forall x\in \R^n,
\end{equation}
where $\lambda\in (0,1)$ and $F^l,F^u:\R^n\to \overline{\R}$ are proper convex and lower semicontinuous extended real valued functions. We call the function $\phi_\lambda$ a scalarized function associated with the PCL function $F$. Note that, \(\lambda \) is strictly restricted to the open interval \((0,1)\) so that the $F^l$ and $F^u$ maintains a non-zero influence on the scalarized function $\phi_\lambda$.
\begin{definition}
     Let $x^*\in \R^n$. Then $x^*$ is referred to an \emph{optimal solution} of (SOP), if there does not exist any other point $x\in \R^n$, such that
    \begin{equation*}
        \phi_\lambda(x)\le \phi_\lambda(x^*).
    \end{equation*}
\end{definition}

\noindent The following results establish the relationship between the two problems (IVOP) and (SOP).
\begin{proposition}\label{ppsn1}
    If $x^*\in\R^n$ is an optimal solution of \textnormal{(SOP)}, then $x^*$ is an $LU$-nondominated solution of \textnormal{(IVOP)}. 
\end{proposition}
\begin{proof}
    Suppose on contrary that $x^*$ is not an $LU$-nondominated solution of (IVOP), then there exists a point $x\in\R^n$, such that
    \begin{equation*}
        F(x)\prec_{LU}F(x^*),
    \end{equation*}
    equivalently,
\[
   \begin{cases}
    F^l(x)< F^l(x^*)\\
    F^u(x)=F^u(x^*)
\end{cases}
\ , \quad \text{or} \quad
\begin{cases}
    F^l(x)= F^l(x^*)\\
    F^u(x)<F^u(x^*)
\end{cases}
\ , \quad \text{or} \quad\
\begin{cases}
    F^l(x)< F^l(x^*),\\
    F^u(x)<F^u(x^*)
\end{cases}
\ .
\]
Now since $\lambda\in( 0,1)$, thus, it follows from the above inequalities
\[
\phi_\lambda(x):=\lambda F^l(x)+(1-\lambda)F^u(x)< \lambda F^l(x^*)+(1-\lambda)F^u(x^*):=\phi_\lambda(x^*),
\]
which contradicts the fact that $x^*$ is an optimal solution of (SOP). Thus, $x^*$ is an $LU$-nondominated solution of (IVOP).
\end{proof}

\begin{example}
     Let the interval-valued optimization problem (IVOP1) be defined as
    \begin{equation*}
    \textbf{(IVOP1)} \qquad \min_{x\in \R} F(x):= [1,1]\odot (x-1)^2\oplus [0,2]\odot |x| \oplus [1,3].\label{IVOP}
\end{equation*}
    Thus, the associated scalarized optimization problem is given by
    \begin{align*}
    \textbf{(SOP1)}\qquad \min_{x\in\R}\phi_\lambda(x):=(x-1)^2+2(1-\lambda)|x|+3-2\lambda,\quad \lambda\in(0,1).
    \end{align*}
    Clearly, $\lambda$ is the optimal solution to the problem (SOP1), $\forall \lambda\in (0,1)$.\\
   Now, to prove that $\lambda\in (0,1)\subseteq \R$ is an $LU$-nondominated solution to the problem (IVOP1), we proceed by contrary and assume that there exists $x\ (\neq \lambda)\in \R$ such that 
   \[F(x)\prec_{LU}F(\lambda),\]
   equivalently,
 \[
   \begin{cases}
    (x-1)^2< (\lambda-1)^2\\
    (x-1)^2+2|x|\\
    \quad \quad =(\lambda-1)^2+2\lambda
\end{cases}
,\quad \text{or} \quad
\begin{cases}
    (x-1)^2= (\lambda-1)^2\\
    (x-1)^2+2|x|\\
    \quad \quad < (\lambda-1)^2+2\lambda
\end{cases}
, \quad \text{or}\\
\begin{cases}
    (x-1)^2< (\lambda-1)^2\\
    (x-1)^2+2|x|\\
    \quad \quad <(\lambda-1)^2+2\lambda
\end{cases}
,
\]
which on simplification (in all the three cases), yields $x=\lambda,\ \forall \lambda\in (0,1)$. \\
Thus, there exists no $x\ (\neq\lambda)\in \R$ such that $F(x)\prec_{LU}F(\lambda)$. Hence, $\lambda$ is an $LU$-nondominated solution of the (IVOP1).
\end{example}

The converse of the Proposition~\ref{ppsn1} does not hold in general. The following example provides that every $LU$-nondominated solution of the (IVOP) is not an optimal solution of (SOP).
\begin{example}
    Let $n=1$, and $F:\R\to \overline{\IR}$ be the interval valued objective of (IVOP), and is defined as
    \begin{align*}
        F(x)=\begin{cases}
            [-x,x^2],\quad \ \ x\in[0,1]\\
            [+\infty,+\infty],\ x\notin[0,1].
        \end{cases}
    \end{align*}
    It is easy to verify that the IVF $F$ is a valid PCL function and $x^*=0$ is an $LU$-nondominated solution of the IVOP. Now, the associated scalarized objective of the (SOP) will be
    \[
    \phi_\lambda(x)=\begin{cases}
            -\lambda x+(1-\lambda)x^2,\quad \ \ x\in[0,1]\\
            \infty,\quad \qquad \qquad \qquad \quad  x\notin[0,1].
        \end{cases}
    \]
    Therefore, the optimal solution of the (SOP) is \[x^*_\lambda=\min\left\{\frac{\lambda}{2(1-\lambda)},\ 1\right \}\neq0,
    \]
    for every $\lambda\in(0,1)$. Consequently, the $LU$-nondominated solution $x^*_\lambda=0$ of (IVOP) is not an optimal solution of (SOP).
\end{example}

\begin{remark}\label{Remark 3.1}
    Note that the Proposition~\ref{ppsn1} depends on the optimal solution of (SOP) to find the $LU$-nondominated solution of (IVOP) and does not rely on the property of convexity, properness and $gH$-lower semicontinuity of the objective function of (IVOP).
\end{remark}

\begin{proposition}\label{eff dom ppsn}
    Let $F:\R^n\to \overline{\IR}$ be an extended real interval valued function and let $\phi_\lambda:\R^n \to \overline{\R}$ be the associated scalarized function, given by $\phi_\lambda(x)=\lambda F^l(x)+(1-\lambda)F^u(x)$ for all $x\in\R^n$ and $\lambda\in(0,1)$. Then 
    \[
    dom\ F=dom\ \phi_\lambda.
    \]
\end{proposition}
\begin{proof}
    Let $x\in dom\ F$, then by the Definition~\ref{effective domain ivf def}, we have $x\in dom\ F^l\cap dom\ F^u$, and therefore
    \[
    \phi_\lambda(x)=\lambda F^l(x)+(1-\lambda)F^u(x)<\infty,
    \]
    which means $x\in dom\ \phi_\lambda$, and consequently $$dom\ F\subseteq dom\ \phi_\lambda.$$ \\
    Conversely, let $x\in dom\ \phi_\lambda$, then
    \[
    \phi_\lambda(x)<\infty
    \]
    \[
    \text{or,}\qquad \lambda F^l(x)+(1-\lambda)F^u(x)<\infty.
    \]
    Since $\lambda\in(0,1)$, it follows from the above inequality
    \begin{align*}
    &F^l(x)<\infty\quad \text{and}\quad F^u(x)<\infty\\
    \text{or,}\qquad &|F^l(x)|<\infty\quad \text{and}\quad |F^u(x)|<\infty\\
    \text{or,}\qquad &\|F(x)\|_{\IR}=\max \{|F^l(x)|,\ |F^u(x)|\}<\infty,
    \end{align*}
     which by Definition~\ref{effective domain ivf def} yields $x\in dom\ F$, and therefore
     $$ dom\ \phi_\lambda\subseteq dom\ F.$$
     Hence, $dom\ F=dom\ \phi_\lambda$.
\end{proof}

\begin{proposition}\label{ppsn:2}
    Let $F:\R^n\to \overline{\IR}$ be an extended real interval valued function and let $\phi_\lambda:\R^n \to \overline{\R}$ be the associated scalarized function, given by $\phi_\lambda(x)=\lambda F^l(x)+(1-\lambda)F^u(x)$ for all $x\in\R^n$ and $\lambda\in(0,1)$, then 
    \begin{enumerate}[label=(\roman*)]
        \item If $F$ is convex $\Rightarrow$ $\phi_\lambda$ is a convex. 
        \item If $F$ is strongly convex with modulus $\gamma>0$ $\Rightarrow$ $\phi_\lambda$ is $\gamma$-strongly convex.
        \item If $F$ is proper $\Rightarrow$ $\phi_\lambda$ is a proper.
        \item If $F$ is $gH$-lower semicontinuous $\Rightarrow$ $\phi_\lambda$ is lower semicontinuous.
    \end{enumerate}
\end{proposition}
\begin{proof}
   (i) Let $F$ be a convex IVF, then by Remark \ref{remark convex}, the following inequalities
       \begin{align*}
           F^l(\theta x+(1-\theta)y)&\le\theta F^l(x)+ (1-\theta) F^l(y),\\
           F^u(\theta x+(1-\theta)y)&\le\theta F^u(x)+ (1-\theta) F^u(y),
       \end{align*}
   holds $\forall x,y\in dom\ F$ and $\theta\in[0,1]$.\\
   Now, since $\lambda\in(0,1)$, on multiplying by $\lambda$ and $(1-\lambda)$ in the above inequalities, respectively, and then summing, we get
    \begin{equation*}
    \begin{aligned}
     \lambda F^l(\theta x+(1-\theta)y)&+(1-\lambda)F^u(\theta x+(1-\theta)y)\\
     &\le \theta \lambda F^l(x)+ (1-\theta)\lambda F^l(y)+\theta (1-\lambda)F^u(x)+ (1-\theta)(1-\lambda)F^u(y),
     \end{aligned}
    \end{equation*}
     equivalently,
    \[
    \phi_\lambda(\theta x+(1-\theta)y)\le \theta \phi_\lambda (x) +(1-\theta)\phi_\lambda (y), \  \forall x,y \in dom\ \phi_\lambda, \ (\text{by Proposition}~\ref{eff dom ppsn}),
    \]
    which shows that $\phi_\lambda$ is convex on $\R^n$.\\
The proof of other parts follows along the similar lines of part (i) as per the assumptions given therein, respectively, and hence being omitted.
\end{proof}

\begin{proposition}\label{gh subdiff ppsn-1}
    Let $F:\R^n\to \overline{\IR}$ be a PCL function and let $\phi_\lambda:\R^n \to \overline{\R}$ be the associated scalarized function as defined in (\ref{scalar function}). Also, let $\hat{\mathcal G}\in \partial_{gH}F(x^*)$ for all $x^*\in int(dom\ F)$, then 
    \[
     \lambda g^-+(1-\lambda)g^+\in \partial \phi_\lambda(x^*),\quad \forall x^*\in int(dom\ \phi_\lambda),
    \]
    where $\hat{\mathcal G}=[g^-,g^+]\in \IR^n$ with $g^- =(g^- _1,g^- _2,\dots,g^- _n),\ g^+ =(g^+ _1,g^+ _2,\dots,g^+ _n)\in \R^n$.
    
\end{proposition}
\begin{proof}
Since $\hat{\mathcal G}\in \partial_{gH}F(x^*)$, then by Remark.~\ref{gh-subdiff remark}, one has
\begin{align*}
 g^-\in\partial F^L(x^*)\ \text{and}\ g^+\in\partial F^U(x^*),
\end{align*}
which by Definition~\ref{def: subdiff}, yields
\begin{align*}
d^\top g^-&\le F^l(x^*+d)-F^l(x^*),\\
d^\top g^+&\le F^u(x^*+d)-F^u(x^*),
\end{align*}
for every $d\in \R^n$. On multiplying by $\lambda$ and $(1-\lambda)$ in the above inequalities, respectively, and then summing, we get
\begin{align*}
  d^\top (\lambda g^-+(1-\lambda)g^+)&\le \lambda F^l(x^*+d)+(1-\lambda)F^u(x^*+d)\\
 &\quad \qquad \qquad \quad \quad  -\lambda F^l(x^*)-(1-\lambda)F^u(x^*),\\
&= \phi_\lambda(x^*+d)-\phi_\lambda(x^*),
\end{align*}
\[
\Rightarrow \ \lambda g^-+(1-\lambda)g^+\in \partial \phi_\lambda(x^*), \quad \forall x^*\in int(dom\ F).
\]
Since, $dom\ F=dom\ \phi_\lambda$, the above inclusion holds for every $x^*\in int(dom\ \phi_\lambda)$.
\end{proof}

\begin{remark}
    The converse of the above result is not generally valid, which means that a point $v\in\partial\phi_\lambda(x^*)$, where $x^*\in int(dom\ \phi_\lambda)$, does not necessarily imply the existence of $\mathcal{\hat G}=[g^-,g^+]\in \partial_{gH}F(x^*)$ such that $v=\lambda g^-+(1-\lambda)g^+$.
\end{remark}

\begin{corollary}Let $x^*\in int(dom\ F)$, and if $\hat 0\in \partial_{gH}F(x^*)$, where $\hat 0=[0,0] \ \text{with}\  0=(0,0,\dots,0)\in \R^n$, then $x^*$ is an $LU$-nondominated solution of the (IVOP). 
\end{corollary}
\begin{proof}
     Since $\hat 0\in \partial_{gH}F(x^*)$, by the above Proposition~\ref{gh subdiff ppsn-1}, one has
     \[
     0 \in \partial\phi_\lambda(x^*). 
     \]
     Since $\phi_\lambda$ is a proper convex function, we have $x^*$ is an optimal solution of the (SOP) (see, for instance, Rockafellar \cite{tyrrell1970convex}), and therefore, by Proposition~\ref{ppsn1}, $x^*$ is an $LU$-nondominated solution of the (IVOP).
\end{proof}

\noindent Based on the above discussion, the (IVOP) can be equivalently studied through its associated (SOP) and it is important to note that the differentiability of the interval-valued objective is not assumed. In particular, the function may fail to be $gH$-differentiable, or its endpoint functions maybe non-smooth, leading to non-applicability of the classical gradient based algorithms. In contrast, the proximal point algorithm is well suitable as it relies only on convexity, properness, and lower semicontinuity, without requiring any differentiability assumptions.
See, for instance Parikh and Boyd \cite{parikh2014proximal}.\\

\noindent We now outline the proximal point algorithm for the (IVOP). The approach is based on the transformation of the the (IVOP) into an associated (SOP). Starting from an initial point $x^0\in dom\ F$, the algorithm computes a sequence of iterates $\{x^k\}$, by solving repeatedly a regularized convex subproblem called the proximal subproblem, such as
\begin{equation}\label{eq:3}
\begin{aligned}
 x^{k+1} &=  \operatorname{prox}_{\mu_k \phi_\lambda}(x^k)\\
        &=\arg\min_{x\in\R^n}\{\phi_\lambda(x)+\frac{1}{2\mu_k}\|x-x^k\|^2\},
\end{aligned}
\end{equation}
where $\mu_k>0$ denotes the proximal parameter satisfying $\sum_{k=0}^{+\infty}\mu_k=+\infty$. The quadratic term enforces proximity to the current iteration and guarantees strong convexity of each subproblem, even when $\phi_\lambda$ is nonsmooth. The complete proximal point algorithm is presented below.

\begin{algorithm}[H]
\caption{Proximal Point Algorithm for (IVOP) \textbf{(PPA-IVOP)}}
\label{alg:1}
\begin{algorithmic}[1]
\Require PCL function $F=[F^l,F^u]$;
         scalarization weight $\lambda\in(0,1)$;
         initial point $x^0\in dom\ F$;
         sequence $\{\mu_k\}_{k\ge 0}\subset[\mu_{\min},\mu_{\max}]$
         with $0<\mu_{\min}\le\mu_{\max}<\infty$;
         tolerance $\varepsilon>0$.
\State Set $k=0$.
\While{$\|x^{k+1}-x^k\|>\varepsilon$}
  \State Compute
         $x^{k+1}=\arg\min_{x\in\R^n}\{\phi_\lambda(x)+\frac{1}{2\mu_k}\|x-x^k\|^2\}$.
  \State Set $k\gets k+1$.
\EndWhile
\State \Return $x^*=x^k$ as an approximate $LU$-nondominated solution.
\end{algorithmic}
\end{algorithm}

\begin{remark} If the endpoints of the PCL function $F$ are differentiable, then the
proximal step ($3$) of the PPA-IVOP reduces to implicit gradient step
$x^{k+1}=x^k-\mu_k\nabla\phi_\lambda(x^{k+1})$. In particular, the proximal update Eq.~\eqref{eq:3} 
satisfies the first order optimality condition \cite{parikh2014proximal},
\[
\nabla \phi_\lambda(x^{k+1})+\frac{1}{\mu_k}(x^{k+1}-x^k)=0.
\]
Furthermore, when the PCL function $F$ is $gH$-differentiable, we cannot guarantee this gradient step, as the $gH$ -differentiability of the IVF $F$ generally does not imply the differentiability of the endpoints of $F$ (See, for instance Chalco Cano $et.al$ \cite{Chalco2011}).

\end{remark}

\section{Convergence Analysis}\label{Section Convergence} 
This section focuses on the convergence behavior of the sequence $\{x^k\}$ generated by PPA-IVOP, and the order of convergence under the following standard assumptions. \\

\noindent \textbf{A1:} The IVF $F=[F^l,F^u]$ is a PCL function with $dom\ F=\R^n$, $i.e.,$ the endpoint functions are finite everywhere.\\
\noindent \textbf{A2:} The scalarization parameter $\lambda\in(0,1)$ is fixed for the statement of each theorem.\\
\textbf{A3:} $\mathcal{X}^{*}_{\lambda}=\arg \min_{x\in \R^n}\phi_{\lambda}(x)\neq \emptyset,$ where $\phi_\lambda$ is defined as in \eqref{scalar function}.\\
\noindent \textbf{A4:} The proximal parameter $\{\mu_k\}\subset [\mu_{\min},\mu_{\max}]$ with $0<\mu_{\min}\le \mu_{\max}<\infty$ is as in PPA-IVOP. In particular, for the rate results, we take $\mu_k\equiv\mu>0$ constant.

\begin{lemma}\label{lem:existence_uniqueness} 
    For each $k=0,1,\dots,K$ where $K\in \mathbb{N}$, the proximal subproblem~(\ref{eq:3}) is well defined.
\end{lemma}

\begin{proof}
Since the IVF $F=[F^l,F^u]$ is a PCL function, the associated scalar function $\phi_\lambda$, defined in (\ref{scalar function}), is proper convex and lower semicontinuous (refer to Proposition \ref{ppsn:2}), which in addition with the regularized quadratic term $\frac{1}{2\mu_k}\|x-x^k\|^2$, yields the proximal subproblem~(\ref{eq:3}) is strongly convex, and therefore it is well defined and has a unique minimum for every $k=0,1,\dots,K$.
\end{proof}

\begin{lemma}\label{lemma endpoint Lipschtz}
    Under \textbf{\textup{(A1)}}, the endpoint functions of the IVF $F$ are locally Lipschitz continuous on $\R^n$. In particular, for every bounded set $\mathcal B\subset \R^n$, there exists $L_{\mathcal B}>0$ such that 
    \begin{align*}
    |F^l(x)-F^l(y)|\le L_{\mathcal B}\|x-y\|\quad \text{and} \quad |F^u(x)-F^u(y)|\le L_{\mathcal B}\|x-y\|,\ \forall x,y\in \mathcal B.
    \end{align*}
\end{lemma}
\begin{proof}
     Since the IVF $F$ is a PCL function, from Proposition~\ref{ppsn:2}, $F^l$ and $F^u$ are convex functions. Therefore, under (\textbf{A1}), the endpoint functions are locally Lipschitz on $\R^n$, (see Theorem 10.4, Convex Analysis, Rockafellar~\cite{tyrrell1970convex}) which means for a bounded set $\mathcal B\subset \R^n$, there exists $L_1,\ L_2>0$, such that 
    \begin{align*}
    |F^l(x)-F^l(y)|\le L_{1}\|x-y\|\ \text{and} \ |F^u(x)-F^u(y)|\le L_{2}\|x-y\|,\ \forall x,y\in \mathcal B.
    \end{align*}
    If we choose $L_{\mathcal B}=\max\{L_1, L_2\}$, we have
    \[
    |F^l(x)-F^l(y)|\le L_{\mathcal B}\|x-y\|\ \text{and} \ |F^u(x)-F^u(y)|\le L_{\mathcal B}\|x-y\|,\ \forall x,y\in \mathcal B.
    \]
\end{proof}
\begin{corollary}\label{corollary IVF F}
    If a sequence $\{x^k\}\subset \R^n$ stays inside a bounded set $\mathcal B$, then for every $x^*\in \mathcal B$,
    \[
    \|F(x^k)\ominus_{gH}F(x^*)\|_{\IR}\le L_{\mathcal B}\|x^k-x^*\|,
    \]
    for every $k=1,2,\dots,K$, where $K\in \mathbb{N}$.
\end{corollary}
\begin{proof}
    The proof follows from the Definitions~\ref{IVF norm} and the above stated Lemma~\ref{lemma endpoint Lipschtz}.
\end{proof}

\begin{theorem}
\label{thm:global_conv}
Let $F:\R^n\to \overline{\IR}$ be a PCL function, and $\phi_\lambda:\R^n \to \overline{\R}$ be the associated scalarized function as defined in (\ref{scalar function}). Let $\mu_k>0$ be the proximal parameter satisfying
$\sum_{k=0}^\infty\mu_k=+\infty$. Then under \textbf{\textup{(A1)--(A3)}}, the sequence $\{x^k\}$ generated by PPA-IVOP satisfies:
\begin{enumerate}[label=(\roman*)]
\item $\phi_{\lambda}(x^{k+1})\le \phi_{\lambda}(x^k)$, for $k\ge 0.$\label{thm glob conv part i}
\item (Fej\'er monotonicity) For every $x^*\in \mathcal X_{\lambda}^*$ and $k\ge 0$,
\begin{align*}
\|x^{k+1}-x^*\|^2&\le \|x^k-x^*\|^2-\|x^{k+1}-x^k\|^2,\\
&\le \|x^k-x^*\|^2,
\end{align*}
and consequently, $\sum_{k=0}^\infty\|x^{k+1}-x^k\|^2\le \|x^0-x^*\|^2<\infty.$\label{thm glob part ii}
\item $\{x^k\}$ converges to a cluster point $\overline{x}_\lambda\in \mathcal X_\lambda^*$.\label{thm glob part iii}
\item $\|F(x^k)\ominus_{gH}F(\overline{x}_\lambda)\|_{\IR}\to 0$ as $k\to \infty$.\label{thm glob part iv}
\item $\overline{x}_\lambda$ is an $LU$-nondominated solution of (IVOP).\label{thm glob part v}
\end{enumerate}
\end{theorem}
\begin{proof}
\begin{enumerate}[label=(\roman*)]
\item[(\textit{i})]By the definition of the proximal subproblem~\eqref{eq:3}, $x^{k+1}$ minimizes
$\phi_\lambda(x)+\frac{1}{2\mu_k}\|x-x^k\|^2$. Substituting $x=x^k$, we have
\begin{align*}
  \phi_\lambda(x^{k+1})+\frac{1}{2\mu_k}\|x^{k+1}-x^k\|^2
  &\le\phi_\lambda(x^k)+\frac{1}{2\mu_k}\|x^{k}-x^k\|^2\\
  &=\phi_\lambda(x^k).
\end{align*}
Since the term $\frac{1}{2\mu_k}\|x^{k+1}-x^k\|^2$ is non-negative for every $k=0,1,\dots,K$, we have
\[
\phi_\lambda(x^{k+1})\le\phi_\lambda(x^k).
\]

\item[(\textit{ii})] From the first order optimality conditions of the proximal subproblem~\eqref{eq:3}, one has $g^{k+1}\in \partial \phi_\lambda(x^{k+1})$, and therefore
\begin{align*}
    &g^{k+1}+\frac{1}{\mu_k}(x^{k+1}-x^k)=0,\\
    \text{or, } \quad &x^k-x^{k+1}=\mu_{k}g^{k+1}, \ \forall k\ge 0. \tag{3}\label{subgrad eq}
\end{align*}
Also, since $x^*\in \mathcal X_\lambda^*$ one has $0\in \partial \phi_\lambda(x^*)$, and therefore from the monotonicity of the subdifferentials, we have
\begin{align*}
\langle g^{k+1},x^{k+1}-x^*\rangle\ge 0.\tag{4}\label{sign eq}
\end{align*}
On expanding $\|x^k-x^*\|^2$ and using~\eqref{subgrad eq} and~\eqref{sign eq}, we obtain
\begin{align*}
    \|x^k-x^*\|^2&=\|x^k-x^{k+1}+x^{k+1}-x^*\|^2,\\
    &=\|x^k-x^{k+1}\|^2+2\mu_k\langle g^{k+1},x^{k+1}-x^*\rangle+\|x^{k}-x^{k+1}\|^2,\\
    &\ge \|x^k-x^{k+1}\|^2+\|x^{k+1}-x^*\|^2,
\end{align*}
equivalently,
\begin{align*}
\|x^{k+1}-x^*\|^2&\le \|x^k-x^*\|^2-\|x^k-x^{k+1}\|^2,\\
&\le \|x^k-x^*\|^2.
\end{align*}
Now, telescoping over $k=0$ to $K\in \mathbb{N}$, we have
\begin{align*}
\sum_{k=0}^K\|x^{k+1}-x^k\|^2&\le\sum_{k=0}^K \|x^k-x^*\|^2-\sum_{k=0}^K\|x^{k+1}-x^*\|^2,\\
&=\|x^0-x^*\|^2-\|x^{K}-x^*\|^2.
\end{align*}
Further, on taking $K\to \infty$, we obtain
\begin{align*}
    \sum_{k=0}^\infty\|x^{k+1}-x^k\|^2\le \|x^0-x^*\|^2<\infty.
\end{align*}

\item [(\textit{iii})] By part~\ref{thm glob part ii}, the sequence $\{x^k\}$ is Fej\'er monotone with respect to $\mathcal X_\lambda^*$ and thus, it follows that $\|x^{k+1}-x^k\|\to 0$, and therefore, from equation~\eqref{subgrad eq}, we have
\begin{align*}\label{eq-1}
\|g^{k+1}\|=\frac{\|x^k-x^{k+1}\|}{\mu_k}\le\frac{\|x^k-x^{k+1}\|}{\mu_{\min}}\to 0,\ \text{as } k\to \infty.\tag{5}
\end{align*}
Now, let $\overline{x}_\lambda\in \R^n$ be any cluster point of the sequence $\{x^k\}$, then for any subsequence $\{x^{k_j+1}\}$ of $\{x^k\}$, one has $x^{k_j+1}\to \overline{x}_\lambda$, and therefore, 
\begin{align*}
    |\phi_\lambda(x^{k_j+1})-\phi_\lambda(\overline{x}_\lambda)|
    &=\lambda|F^l(x^{k_j+1}-F^l(\overline{x}_\lambda)|+(1-\lambda)|F^u(x^{k_j+1}-F^u(\overline{x}_\lambda)|,\\
    &\le L_\mathcal B\|x^{k_j+1}-\overline{x}_\lambda\|\to 0\ \text{as } j\to \infty.\ (by\ using\ \text{Lemma}~\ref{lemma endpoint Lipschtz})\\
    \Rightarrow \phi_\lambda(x^{k_j+1})&\to\phi_\lambda(\overline{x}_\lambda),\ \text{as } j\to \infty. \tag{6}\label{eq-2}
    \end{align*}
On using the subgradient of $\phi_\lambda$ at $x^{k_j+1}\in \R^n$, we obtain
\begin{align*}
\phi_\lambda(y)&\ge \phi_\lambda(x^{k_j+1})+\langle g^{k_j+1},y-x^{k_j+1}\rangle,\ \forall y\in \R^n,
\end{align*}
which by the implications \eqref{eq-1} and \eqref{eq-2} as $j\to \infty$, yields
\begin{align*}
&\phi_\lambda(y)\ge \phi_\lambda(\overline{x}_\lambda),\ \forall y\in \R^n.\\
\Rightarrow \ &\overline{x}_\lambda\in \mathcal X_\lambda^*,
\end{align*}
and we conclude that the sequence $\{x^k\}$  converges to $\overline{x}_\lambda\in \mathcal X_\lambda^*$. 

\item[(\textit{iv})] By part~\ref{thm glob part iii}, the sequence $\{x^k\}$ is convergent and hence bounded, therefore there exists a constant $\mathcal R>0$ such that
\[
\|x^k-\overline{x}_\lambda\|< \mathcal R,\ \forall k\ge 0 \text{ and } \overline{x}_\lambda\in \mathcal X_\lambda^*.
\]
Let $\mathcal B=\{x\in \R^n| \ \|x-\overline{x}_\lambda\|< \mathcal R\}$ be a bounded set containing $\overline{x}_\lambda$ and every iterate $x^k$, then by Corollary~\ref{corollary IVF F} there exist a constant $L_{\mathcal B}>0$, such that
\[
\|F(x^k)\ominus_{gH}F(\overline{x}_\lambda)\|_{\IR}\le L_{\mathcal B}\|x^k-\overline{x}_\lambda\|\to 0 \text{ as }k\to \infty.
\]
It follows that
\[
\|F(x^k)\ominus_{gH}F(\overline{x}_\lambda)\|_{\IR}\to 0\ \text{as}\ k\to\infty.
\]

\item[(\textit{v})] Since $\overline{x}_\lambda\in \mathcal X_\lambda^*$ is an optimal solution of (SOP), therefore from Proposition~\ref{ppsn1}, $\overline{x}_\lambda$ is an $LU$-nondominated solution of (IVOP).
\end{enumerate}
\end{proof}

\begin{remark}
    The above theorem holds for any initial point $x^0\in \R^n$ and any $\lambda\in (0,1)$, so different values of $\lambda$ generally yields different $LU$-nondominated solutions. By iterating PPA-IVOP over a grid $\{\lambda_i\}_{i=1}^m$ one obtain a collection of $\{x^*(\lambda_i)\}_{i=1}^m$ of $LU$-nondominated solutions. However, under an additional structural assumption on the interval valued objective of (IVOP), the PPA-IVOP generates unique minimizer, independent of the choice of $\lambda\in (0,1).$
\end{remark}
\begin{proposition}\label{convg structure ppsn}
Let $F:\R^n\to \overline{\IR}$ be a PCL function and if its endpoint functions $F^l$ and $F^u$ possesses a unique common minimizer, say $x^*\in \R^n$. Then $x^*$ is the unique $LU$-nondominated solution of the IVOP and consequently, the sequence $\{x^k\}$ generated by PPA-IVOP converges to $x^*$, independent of the choice of $\lambda\in (0,1)$.
\end{proposition}
\begin{proof}
   Since $x^*$ is a unique minimizer of the endpoint functions $F^l$ and $F^u$, one has
    \[
     \arg \min_{x\in \R^n}F^l(x)=\arg \min_{X\in \R^n}F^u(x)=\{x^*\},
    \]
    \[
    \text{or,}\ F ^l(x^*)< F^l(x)\quad \text{and}\quad F^u(x^*)< F^u(x), \ \forall x\ (\neq x^*)\in \R^n.
    \]
    On multiplying by $\lambda$ and $(1-\lambda)$ and, then summing, we get
    \[
    \phi_\lambda(x^*)< \phi_\lambda(x),\ \forall x\ (\neq x^*) \in \R ^n,\ \forall \lambda\in(0,1),
    \]
    and therefore, one has\[
    \arg \min_{x\in\R^n}\phi_\lambda(x)=\{x^*\}, \ \forall \lambda\in(0,1),
    \]
   which shows that $x^*$ is a unique optimal solution of the (SOP), and consequently, a unique $LU$-nondominated solution of the (IVOP) by Proposition~\ref{ppsn1}. Further, by Theorem,~\ref{thm:global_conv}, since every cluster point of the sequence $\{x^k\}$ generated by PPA-IVOP converges to the optimal solution of the (SOP), and hence the sequence $\{x^k\}$ converges to a unique minimizer $x^*$, independent of the choice of $\lambda\in(0,1)$.  
\end{proof}

We now establish the rate of convergence of PPA-IVOP for solving the problem (IVOP).

\begin{theorem}\label{thm ergodic rate}
    Assume that the assumptions \textbf{\textup{(A1)--(A4)}} and the hypothesis of Theorem~\ref{thm:global_conv} holds. Then for every $x^*\in \mathcal X_\lambda^*$ and $K\in \mathbb N$, we have
    \[
    \min_{0\le k\le K-1}\big (\phi_\lambda(x^{k+1})-\phi_\lambda(x^*)\big )\le \frac{\|x^0-x^*\|^2}{2\mu K}.
    \]
\end{theorem}
\begin{proof}
From the subdifferentiablity of the function $\phi_\lambda$ at the point $x^{k+1}\in \R^n$, we have
    \begin{align*}
        \phi_\lambda(x^*)\ge \phi_\lambda(x^{k+1})+\langle g^{k+1},x^*-x^{k+1}\rangle, \ \forall k\ge 0,\ x^*\in \mathcal X_\lambda^*,\ \forall g^{k+1}\in \partial \phi_\lambda(x^{k+1}),
    \end{align*}
    which together with the equation \eqref{subgrad eq}, yields
    \begin{align*}
    \phi_\lambda(x^{k+1})-\phi_\lambda(x^*)&\le \frac{1}{\mu}\langle x^k-x^{k+1}, x^{k+1}-x^*\rangle,\\
    &=\frac{1}{2\mu}\big (\|x^k-x^*\|^2-\|x^k-x^{k+1}\|^2-\|x^{k+1}-x^*\|^2\big ),\\
    &\le \frac{1}{2\mu}\big (\|x^k-x^*\|^2-\|x^{k+1}-x^*\|^2\big ).\ (since\ \|x^k-x^{k+1}\|^2\ge 0)
    \end{align*}
On telescoping over $k=0,\dots,K-1$, we obtain
\begin{align*}
    \sum_{k=0}^{K-1}\big ( \phi_\lambda(x^{k+1})-\phi_\lambda(x^*)\big )&\le \frac{1}{2\mu}\sum_{k=0}^{K-1}\big (\|x^k-x^*\|^2-\|x^{k+1}-x^*\|^2\big ),\\
    &=\frac{1}{2\mu}\big (\|x^0-x^*\|^2-\|x^{K}-x^*\|^2\big ),\\
    &\le \frac{1}{2\mu}\|x^0-x^*\|^2.\ (since\ \|x^K-x^*\|^2\ge 0)
\end{align*}
By Theorem~\ref{thm:global_conv}\ref{thm glob conv part i}, one has the gap $\phi_\lambda(x^{k+1})-\phi_\lambda(x^*)$ is non-increasing. Consequently,
\begin{align*}
K. \min_{0\le k\le K-1} \big (\phi_\lambda(x^{k+1})-\phi_\lambda(x^*)\big )&\le \sum_{k=0}^{K-1}\big ( \phi_\lambda(x^{k+1})-\phi_\lambda(x^*)\big ),\\
&\le \frac{1}{2\mu}\big (\|x^0-x^*\|^2\big ),\\
\text{or, }\min_{0\le k\le K-1} \big (\phi_\lambda(x^{k+1})-\phi_\lambda(x^*)\big )&\le \frac{1}{2\mu K}\big (\|x^0-x^*\|^2\big ).
\end{align*}
\end{proof}
\begin{corollary}
    In addition to assumptions made in Theorem~\ref{thm ergodic rate}, if the Proposition~\ref{convg structure ppsn} also holds, then for every $\lambda \in (0,1)$ and $K\in \mathbb{N}$,\[
    \min_{0\le k\le K-1}\|F(x^{k+1})\ominus_{gH}F(x^*)\|_\IR\le \frac{C(\lambda)}{2\mu K}\|x^0-x^*\|^2,
    \]
    where $C(\lambda):=\max \big\{\frac {1}{\lambda},\frac {1}{1-\lambda}\big \}$ and $x^*\in \mathcal X_\lambda^*$.
\end{corollary}
\begin{proof}
Let $\{x^k\}$ be the sequence generated by the PPA-IVOP, then from Proposition~\ref{convg structure ppsn}, we have
    \begin{align*}
        F^l(x^{k+1})-F^l(x^*)\ge 0\quad \text{and}\quad F^u(x^{k+1})-F^u(x^*)\ge 0,\ \forall k\ge 0,
    \end{align*}
which on scalarization in terms of $\phi_\lambda$, yields
\begin{align*}
    \phi_\lambda(x^{k+1})-\phi_\lambda(x^*)=\lambda\big ( F^l(x^{k+1})-F^l(x^*)\big )+(1-\lambda)\big ( F^u(x^{k+1})-F^u(x^*)\big ),
\end{align*}
    and consequently,
    \begin{align*}
    0\le F^l(x^{k+1})-F^l(x^*)&\le \frac{1}{\lambda}\big(\phi_\lambda(x^{k+1})-\phi_\lambda(x^*)\big )\\
        \text{and  }
       0\le  F^u(x^{k+1})-F^u(x^*)&\le \frac{1}{1-\lambda}\big(\phi_\lambda(x^{k+1})-\phi_\lambda(x^*)\big ), 
    \end{align*} 
    for every $\lambda\in (0,1)$. Now, using the Definition \ref{IVF norm} for the IVF $F$, we obtain
\begin{align*}
    \|F(x^{k+1})\ominus_{gH}F(x^*)\|_\IR&=\max\{|F^l(x^{k+1})-F^l(x^*)|, |F^u(x^{k+1})-F^u(x^*)|\},\\
    &\le C(\lambda)\big (\phi_\lambda(x^{k+1})-\phi_\lambda(x^*)\big ),
\end{align*}
where $C(\lambda):=\max \big\{\frac {1}{\lambda},\frac {1}{1-\lambda}\big \}$. 

\noindent Taking the minimum over $k$ from $k=0,\dots,K-1$, and using the results of Theorem~\ref{thm ergodic rate} in the above inequality, one has
\begin{align*}
\min_{0\le k\le K-1}\|F(x^{k+1})\ominus_{gH}F(x^*)\|_\IR\le \frac{C(\lambda)}{2\mu K}\|x^0-x^*\|^2.
\end{align*}
\end{proof}
We now show that for the strongly convex case, the PPA-IVOP has a linear rate of convergence.
\begin{theorem}\label{thm:linear}
Suppose the IVF $F$ is $\gamma$-strongly convex and $\mu_k\equiv\mu>0$ be a constant proximal parameter. Further, we assume that the hypothesis of Theorem~\ref{thm:global_conv} are fulfilled, then
\[
\|x^{k+1}-x^*\|^2 \le \left(\frac{1}{1+\gamma \mu}\right)^{k+1}\|x^0-x^*\|^2,
\]
for every $k=0,1,\dots,K-1$, where $K\in \mathbb{N}$ and $x^*\in \mathcal X_\lambda^*$. Moreover, there exists a constant $L_\mathcal B>0$ such that 
\[
\|F(x^{k+1})\ominus_{gH}F(x^*)\|_\IR\le L_\mathcal B\left(\frac{1}{1+\gamma \mu}\right)^{\frac{k+1}{2}}\|x^0-x^*\|.
\]
\end{theorem}
\begin{proof}
Since the IVF $F$ is strongly convex, from Proposition~\ref{ppsn:2} the scalarized function $\phi_\lambda$  is strongly convex, and  therefore, from Remark~\ref{sugdiff remark}, we have 
\begin{align*}
\phi_\lambda(x^*)\ge \phi_\lambda(x^{k+1})+\langle g^{k+1}, x^*-x^{k+1}\rangle+\frac{\gamma}{2}\|x^{k+1}-x^*\|^2,\\
\text{or, }\phi_\lambda(x^*)-\phi_\lambda(x^{k+1})
+\langle g^{k+1}, x^{k+1}-x^*\rangle\ge \frac{\gamma}{2}\|x^{k+1}-x^*\|^2.
\end{align*}
Further, since $x^*\in \mathcal X_\lambda^*$, the term $(\phi_\lambda(x^*)-\phi_\lambda(x^{k+1}))$ is non-positive, dropping it in the above inequality yields
\begin{align*}
\langle g^{k+1}, x^{k+1}-x^*\rangle\ge \frac{\gamma}{2}\|x^{k+1}-x^*\|^2.\tag{7}\label{strong subgrad eq}
\end{align*}
Expanding $\|x^{k}-x^*\|^2$ and using \eqref{subgrad eq} and \eqref{strong subgrad eq}, we obtain
\begin{align*}
    \|x^{k}-x^*\|^2&=\|x^k-x^{k+1}+x^{k+1}-x^*\|^2,\\
    &=\|x^{k}-x^{k+1}\|^2+2\mu\langle g^{k+1}, x^{k+1}-x^*\rangle+\|x^{k+1}-x^*\|^2,\\
    &\ge (1+\gamma \mu)\|x^{k+1}-x^*\|^2.\quad (since\  \|x^{k}-x^{k+1}\|^2\ge 0)
\end{align*}
Therefore, iterating the above inequality till $k=0$, yields
\begin{align*}
\|x^{k+1}-x^*\|^2\le \left (\frac{1}{1+\gamma \mu}\right )\|x^k-x^*\|^2&\le  \left (\frac{1}{1+\gamma \mu}\right )^2\|x^{k-1}-x^*\|^2\\
&\le \dots \le  \left (\frac{1}{1+\gamma \mu}\right )^{k+1}\|x^0-x^*\|^2.
\end{align*}
Further, since the sequence $\{x^k\}$ converges to $x^*$, it is bounded. Hence, there exists $\mathcal R>0$ such that the trajectory  is contained in an open ball $\mathcal B:=\{x\in \R^n |\  \|x-x^*\|<\mathcal R\}$. Therefore, by Corollary~\ref{corollary IVF F}, there exist a constant $L_\mathcal B>0$, such that 
\begin{align*}
    \|F(x^{k+1})\ominus_{gH}F(x^*)\|_\IR&\le L_\mathcal B\|x^{k+1}-x^*\|,\\
    &\le L_\mathcal B\left(\frac{1}{1+\gamma \mu}\right)^{\frac{k+1}{2}}\|x^0-x^*\|.
\end{align*}
\end{proof}

\subsection*{Numerical Experiments}\label{Section numerical experiments}

In this section, we illustrate the performance of the proposed PPA-IVOP. The experiments demonstrate the convergence behavior of the algorithm and its ability to compute $LU$-nondominated solutions of the (IVOP). 

For the given experiments, $dom\ F=\R^n$, the proximal parameter sequence $\{\mu_k\}$ is taken as a constant value $\mu=0.2$, and the algorithm is initialized from different starting points $x^0\in dom\ F$ as well as the different values of scalarization parameter $\lambda\in(0,1)$. Further, we assume that the iterative procedure is terminated using the stopping rule $\|x^{k+1}-x^k\|<\varepsilon=10^{-6}$ (tolerance).

\begin{example}\label{example 1}
    We study the interval valued optimization problem, defined as follows
    \[
    \textbf{(IVOP2)}\quad \min_{x\in\R}F(x):=
    [3,4]\odot |x|\oplus[1,6]\odot|x-5|\oplus [4,10],
    \]
    and thus, the associated scalarized optimization problem is given by
    \[
    \textbf{(SOP2)}\quad\min_{x\in\R}\phi_\lambda(x):=(4-\lambda)|x|+(6-5\lambda)|x-5|+10-6\lambda.
    \]
where $\lambda\in(0,1)$.
The graphical representation of the IVF $F$ are exhibited in  Figure~\ref{figure 5.1}, where the blue region is the interval region and the black segment on the $x$-axis indicates the $LU$-nondominated solution set $[0,5]$.
\begin{figure}[H]
\centering
\includegraphics[width=0.5\linewidth]{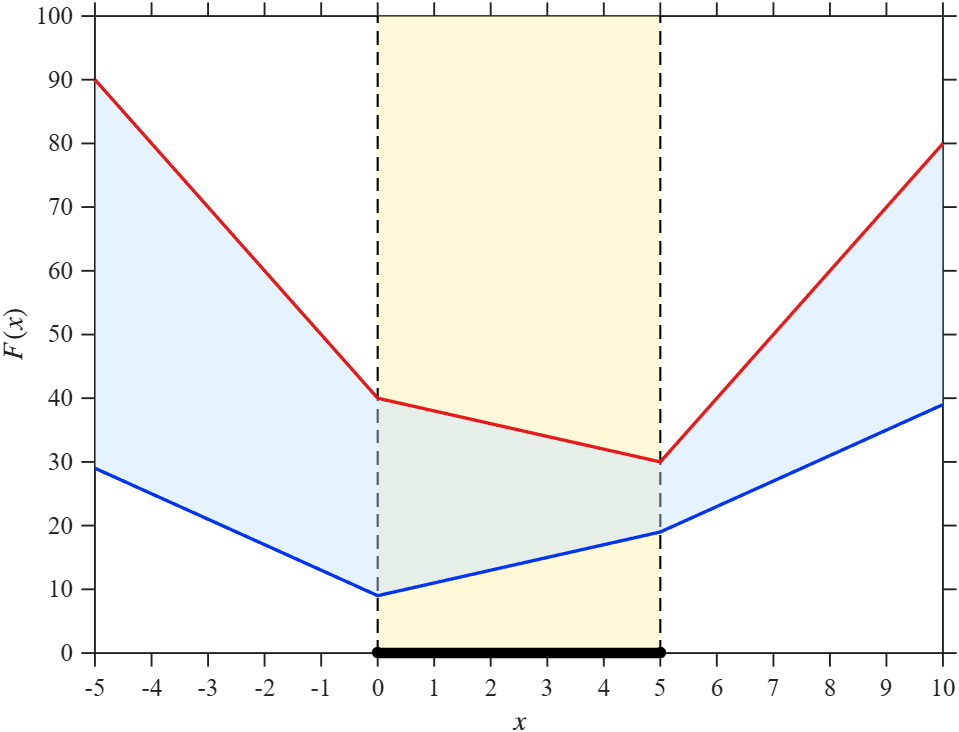}
\caption{Contour plot of the IVF $F$ of (IVOP2).}
\label{figure 5.1}

\end{figure}
 It is easy to verify that the IVF $F$ is a PCL function. The implementation of the PPA-IVOP with different choices of the scalarization parameter $\lambda\in (0,1)$ and arbitrary initial points $x^0\in \R$, is presented in Table~\ref{tab:ex5.1}.

\begin{table}[H]
\centering
\small
\setlength{\tabcolsep}{4pt}
\renewcommand{\arraystretch}{1.15}
\begin{tabular}{cccc}
\toprule
$\lambda$ &
\shortstack{Initial Point \\$x^0$} &
\shortstack{Iterations\\ $(k)$} &
\shortstack{LU-nondominated \\Solution $x^*$} \\
\midrule

\multirow{3}{*}{$0.1$}
& $-4$   & 20 & 5.0000  \\
& $2.5$  & 9 &  5.0000 \\
& $9$    & 4 &  5.0000  \\
\midrule

\multirow{3}{*}{$0.3$}
& $-4$   & 36 &  5.0000\\
& $2.5$  & 17 &  5.0000  \\
& $9$    & 4 &   5.0000 \\
\midrule

\multirow{3}{*}{$0.5$}
& $-4$   & 4 &   0.0011 \\
& $2.5$  & 1 &  2.5000  \\
& $9$    & 4 &  4.9987\\
\midrule

\multirow{3}{*}{$0.7$}
& $-4$   & 5 &  0.0000 \\
& $2.5$  & 17 &  0.0000 \\
& $9$    & 37 &  0.0000  \\
\midrule

\multirow{3}{*}{$0.9$}
& $-4$   & 6 &   0.0000 \\
& $2.5$  & 9 &  0.0000 \\
& $9$    & 22 &  0.0000  \\
\bottomrule
\end{tabular}
\caption{Output of PPA-IVOP for (IVOP2).}
\label{tab:ex5.1}
\end{table}
Since the (IVOP2) admits infinitely many $LU$-nondominated solutions,  Table.~\ref{tab:ex5.1} reports the few of them on applying PPA-IVOP for different values of $\lambda\in(0,1)$. This confirms that the algorithm reliably locates the $LU$-nondominated solutions of (IVOP2), with $\lambda$ acting as the selection mechanism across the nondominated set $[0,5]$.
\end{example}

In the above example, the endpoint functions have different minimizers, thus, producing a set of $LU$-nondominated solutions of the (IVOP2). Consequently, the scalarization parameter $\lambda$ acts as a selection mechanism and the different values of $\lambda$ recovering different $LU$-nondominated solutions contained in that set. On the other hand, in the following example, we consider an (IVOP) such that the endpoint functions of the interval valued objective possesses a common minimizer and therefore satisfies the assumptions of Proposition~\ref{convg structure ppsn}, and thus, the PPA-IVOP converges to a unique $LU$-nondominated solution, independent to any choice of $\lambda\in (0,1)$ and arbitrary initial points $x^0\in dom\  F$.

\begin{example}\label{example 5.2}
    We study the interval valued optimization problem, defined as follows
    \begin{align*}
    \textbf{(IVOP3)}\quad \min_{(x_1,x_2)\in\R^2}F(x_1,x_2):=\mathbf{1}\odot (|x_1-1|+|x_2-2|)\quad \qquad \quad\ \ \ \quad \  \\\oplus [1,2]\odot(x_1-1)^2+(x_2-2)^2\oplus[2,3].\label{IVOP2}
    \end{align*}
    The associated scalarized optimization problem is given by
    \begin{align*}
    \textbf{(SOP3)}\quad \min_{(x_1,x_2)\in \R^2}\phi_\lambda(x_1,x_2):=|x_1-1|&+|x_2-2|\\+&(2-\lambda)\bigl ((x_1-1)^2+(x_2-2)^2\bigr )+(3-\lambda).
    \end{align*}
   It is easy to verify that the IVF $F$ is a PCL function. The implementation of the PPA-IVOP with different choices of the scalarization parameter $\lambda\in(0,1)$ and various initial points $(x_1^0,x_2^0)\in \R^2$, is presented in Table~\ref{tab:ex5.2}.
\begin{table}[H]
\centering
\small
\setlength{\tabcolsep}{4pt}
\renewcommand{\arraystretch}{1.15}
\begin{tabular}{cccc}
\toprule
$\lambda$ &
\shortstack{Initial Point \\$(x_1^0,x_2^0)^\top$} &
\shortstack{Iterations\\ $(k)$} &
\shortstack{LU-nondominated \\Solution $\bar x=(\bar x_1,\bar x_2)^\top$} \\
\midrule

\multirow{3}{*}{$0.1$}
& $(5,5)^\top$   & 6  &  \\
& $(-3,-3)^\top$ & 7  & $(1,2)^\top$ \\
& $(10,-5)^\top$ & 8  &  \\
\midrule

\multirow{3}{*}{$0.3$}
& $(5,5)^\top$   & 7  &  \\
& $(-3,-3)^\top$ & 7  & $(1,2)^\top$ \\
& $(10,-5)^\top$ & 8  &  \\
\midrule

\multirow{3}{*}{$0.5$}
& $(5,5)^\top$   & 7  &  \\
& $(-3,-3)^\top$ & 7  & $(1,2)^\top$ \\
& $(10,-5)^\top$ & 9  &  \\
\midrule

\multirow{3}{*}{$0.7$}
& $(5,5)^\top$   & 7  &  \\
& $(-3,-3)^\top$ & 8  & $(1,2)^\top$ \\
& $(10,-5)^\top$ & 9  &  \\
\midrule

\multirow{3}{*}{$0.9$}
& $(5,5)^\top$   & 8  &  \\
& $(-3,-3)^\top $ & 8  & $(1,2)^\top$ \\
& $(10,-5)^\top$ & 10 &  \\
\bottomrule
\end{tabular}
\caption{Output of PPA-IVOP for (IVOP3).}
\label{tab:ex5.2}
\end{table}

The numerical results reported in Table~\ref{tab:ex5.2} demonstrate that, PPA-IVOP converges to a unique $LU$-nondominated solution $(1,2)^\top$, independent to any choice of $\lambda\in(0,1)$. The number of iterations remain small for all the experiments, illustrating the robustness of the algorithm with respect to both the initial point and the scalarization parameter. These observations are also consistent with the theoretical convergence results established in Theorem~\ref{thm:global_conv}.

\begin{figure}[H]
\centering

\begin{subfigure}{0.48\textwidth}
\centering
\includegraphics[width=0.9\linewidth]{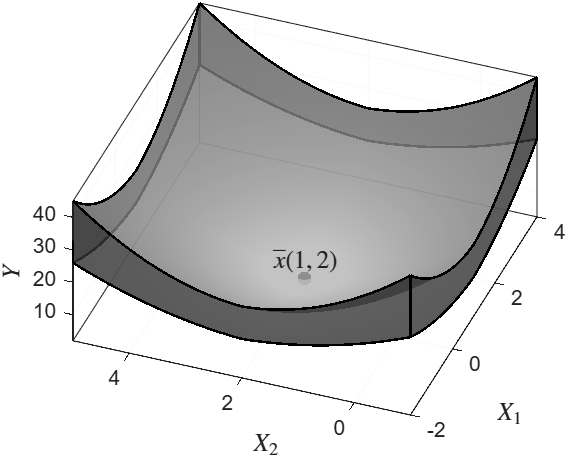}
\caption{Contour plot of the IVF $F$ of (IVOP3).}
\end{subfigure}
\hfill
\begin{subfigure}{0.48\textwidth}
\centering
\includegraphics[width=1.0\linewidth]{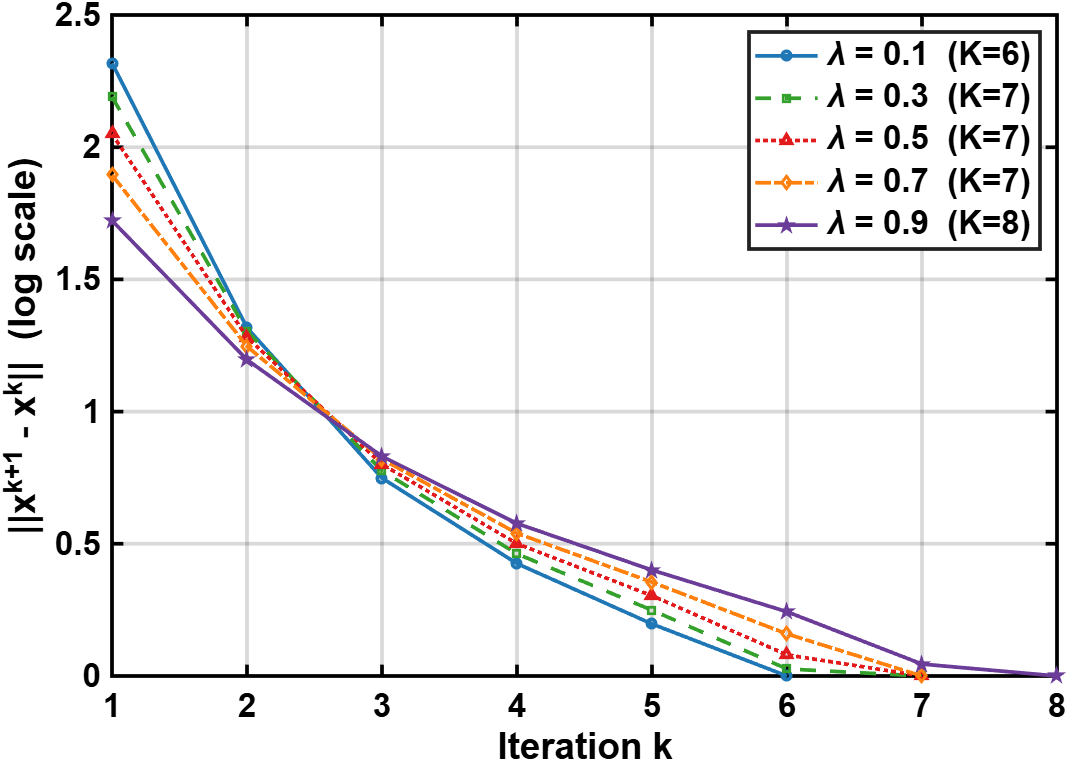}
\caption{Convergence of iterates $(x^k)$.}\label{fig 2b}
\end{subfigure}

\caption{Graphical illustrations of the PPA-IVOP on the (IVOP3) starting from the initial point $x^0=(5,5)^\top$.}
\label{fig:ppa_results}

\end{figure}
\noindent Moreover, one can easily verify that the IVF $F$ is strongly convex, therefore the linear convergence result, established in Theorem~\ref{thm:linear}, is confirmed numerically by the rapid decay of the iterates shown in Figure~\ref{fig 2b}.

\end{example}

\subsection*{Comparative analysis}

To further assess the numerical behavior of the PPA-IVOP, we compare it with the projected $gH$-subgradient method proposed by
Kumar $et.\ al$~\cite{KSR} on their Example~5.1, whose interval form is given by
\[
\min_{x\in[-1,1]}F(x)=\begin{cases}
    [2+|x|,\,6],\quad \quad \quad x\in[-1,1]\\
    [2+|x|,\,4+2|x|]\quad \text{otherwise.}
\end{cases}
\]
The above problem can equivalently be written as an unconstrained (IVOP) given by
\[
\textbf{(IVOP4)}\quad \min_{x\in\R}F(x)=\begin{cases}
    [2+|x|,\,6],\quad \quad \quad x\in[-1,1]\\
    [2+|x|,\,4+2|x|]\ \ x\notin [-1,1].
\end{cases}
\]
The endpoint functions of the IVF $F$ are convex and continuous, so $F$ is a PCL function with $dom\ F=\R$. The unique $LU$-nondominated solution is
$x^*=0$ of the (IVOP4) which coincides with the weak efficient solution of Kumar $et.\ al$~\cite{KSR}, therefore both the
methods target the same point. The proximal parameter $\mu=1$ and tolerance of $10^{-6}$ are selected for the proposed PPA-IVOP, while the parameters of the projected $gH$-subgradient method are taken according to~\cite{KSR}. Tables~\ref{tab:comp ppa} and \ref{tab:comp psm} reports the results obtained by the proposed PPA-IVOP and the projected $gH$-subgradient method on the (IVOP4) problem, respectively.

\begin{table}[H]
\centering
\small
\setlength{\tabcolsep}{4pt}
\renewcommand{\arraystretch}{1.15}
\begin{tabular}{cccc}
\toprule
$\lambda$ &
\shortstack{Initial Point \\$x^0$} &
\shortstack{Iterations\\ $(k)$} &
\shortstack{LU-nondominated \\Solution $x^*$} \\
\midrule

\multirow{3}{*}{$0.1$}
& $-1$   & 11  &  \\
& $0.5$ & 6  & $0$ \\
& $1$ & 11  &  \\
\midrule

\multirow{3}{*}{$0.3$}
& $-1$   & 5  &  \\
& $0.5$ & 3  & $0$ \\
& $1$ & 5  &  \\
\midrule

\multirow{3}{*}{$0.5$}
& $-1$   & 3  &  \\
& $0.5$ & 2  & $0$ \\
& $1$ & 3  &  \\
\midrule

\multirow{3}{*}{$0.7$}
& $-1$   & 3  &  \\
& $0.5$ & 2  & $0$ \\
& $1$ & 3  &  \\
\midrule

\multirow{3}{*}{$0.9$}
& $-1$   & 3  &  \\
& $0.5$ & 2  & $0$ \\
& $1$ & 3  &  \\
\bottomrule
\end{tabular}
\caption{Output of PPA-IVOP for (IVOP4).}
\label{tab:comp ppa}
\end{table}

\begin{table}[H]
\centering
\small
\setlength{\tabcolsep}{4pt}
\renewcommand{\arraystretch}{1.15}
\begin{tabular}{cccc}
\toprule
Step rule $\lambda_k$ &
\shortstack{Initial Point \\$x^0$} &
\shortstack{Iterations\\ $(k)$} &
\shortstack{Solution $x^*$} \\
\midrule

\multirow{3}{*}{$1/(2k)$}
& $-1$   & 18  & $-0.0003$ \\
& $0.5$ & 1  & $0$ \\
& $1$ & 24  &  $0.0030$\\

\bottomrule
\end{tabular}
\caption{Output of projected $gH$-subgradient method for (IVOP4).}
\label{tab:comp psm}
\end{table}
For every $\lambda\in(0,1)$ and every initial point $x^0$, PPA-IVOP converges to the
$LU$-nondominated solution $x^*=0$ exactly, in at most $11$ iterations, and $F(x^*)=[2,6]$ in all runs. Therefore, $\lambda$ acts purely as a
speed parameter which larger value of $\lambda$ giving faster convergence. By comparison, for the step rule $\lambda_k=1/(2k)$, the projected $gH$-subgradient
method requires $18$ and $24$ iterations from $x^0=-1$ and $x^0=1$
and returns only approximations of $x^*$ ($-0.0003$ and $0.0030$), whereas
PPA-IVOP attains $x^*$ exactly in $3$ to $11$ iterations from the same initial points. Moreover, PPA-IVOP has no step-size sequence to tune, and its iterates
are guaranteed to converge (Theorem~\ref{thm:global_conv}), together with a sublinear rate in
the convex case and a linear rate in the strongly convex case
(Theorems~\ref{thm ergodic rate} and~\ref{thm:linear}), whereas for the subgradient method only the convergence of
the objective values are established~\cite{KSR}.

\section{Application in interval valued lasso Problem}\label{Section Application}

In this section, we illustrate the effectiveness of the proposed PPA-IVOP to an interval valued $l_1$-regularized linear regression problem, commonly known as the $lasso$ problem. This extends the classical \textit{lasso} formulation introduced by Tibshirani \cite{Tibshirani1996} to handle uncertain or interval-valued observations.\\

\noindent Consider a dataset consisting of $n$ pairs of interval-valued training dataset $\{X_k,Y_k\}_{k=1}^n$, where for each sample $k = 1, 2, \ldots, n$, the feature vector is represented by the interval $X_k\in\IR^m$ given by $X_k=\{X_k^1,X_k^2,\dots,X_k^m\ | \ X_k^j\in \IR\ ;\  j=1,2,\dots,m\}$ and the corresponding target is $Y_k \in\IR$. The objective is to determine a regression parameter vector $\theta=(\theta_1,\theta_2,\dots,\theta_m)\in \R^m$ that fits the given interval-valued dataset. 

\medskip
\noindent For each interval valued feature, define its center and radius by
\begin{align*}
&X_{k,c}^j=\frac{\underline{X}_k^j+\overline{X}_k^j}{2},\quad X_{k,r}^j=\frac{\overline{X}_k^j-\underline{X}_k^j}{2}\\
\text{and }& Y_{k,c}=\frac{\underline{Y}_k+\overline{Y}_k}{2},\quad Y_{k,r}=\frac{\overline{Y}_k-\underline{Y}_k}{2}.
\end{align*}
For a fixed parameter value $\theta$, the interval valued prediction function $\mathcal P(\cdot;\theta):\IR^m\to\IR$, is defined as
\[
\mathcal P(X_k;\theta)=\bigoplus_{j=1}^{m}(\theta_j\odot X^j_k),
\]
where $\theta_j\odot X^j_k=\left [\theta_j X_{k,c}^j-|\theta_j|X_{k,r}^j,\, \theta_j X_{k,c}^j+|\theta_j|X_{k,r}^j\right]$. Consequently, the prediction function is written as
\[
\mathcal P(X_k;\theta)=\left [\mathcal P_{k,c}(\theta)-\mathcal P_{k,r}(\theta),\,\mathcal P_{k,c}(\theta)+\mathcal P_{k,r}(\theta), \right ],
\]
where $\mathcal P_{k,c}(\theta)=\sum_{j=1}^m\theta_j X_{k,c}^j$ and 
 $\mathcal P_{k,r}(\theta)=\sum_{j=1}^m|\theta_j| X_{k,r}^j$. This function represents the predicted interval corresponding to the $k^{th}$ feature vector $X_k$. The regression parameter $\theta$ is obtained by minimizing an interval valued loss function so that the resulting prediction $\mathcal P(X_k;\theta)$ provides approximation to the target interval $Y_k$, for every $k=1,2,\dots,n$. 
 
  Fixing the training dataset, we define the center-radius based interval valued squared loss function $\mathcal L_k(\theta)$ as
 \[
 \mathcal L_k(\theta)=[(\mathcal P_{k,c}(\theta)-Y_{k,c})^2,\, (\mathcal P_{k,c}(\theta)-Y_{k,c})^2+(\mathcal P_{k,r}(\theta)+Y_{k,r})^2],
 \]
where the lower endpoint measure the squared error between the predicted center and the target center, while the additional term in the upper endpoint accounts for the uncertainty associated with the interval valued data.
 
It is observable that IVF $\mathcal L_k$ is convex and $gH$-continuous. Consequently, this leads to the interval valued lasso optimization problem, defined as follows
\begin{equation*}
\textbf{(IVOP5)}\ \min_{\theta\in\R^m} E(\theta)=\frac{1}{2}\odot\bigoplus_{k=1}^n\mathcal L_k(\theta)\oplus\mathcal{A}\odot\|\theta\|_1,\label{IVOPapp}
\end{equation*}
where $\mathcal A\in\IR_+$ is the interval valued regularization parameter, and $\|.\|_1$ is the $l_1$-norm, defined as
\[
\|\theta\|_1=\sum_{j=1}^m|\theta_j|,\quad \forall j=1,2,\dots,m.
\]
The objective function $E(\theta)$ is a PCL function since both the squared loss IVF $\mathcal L_k$ and the $\ell_1$-norm IVF $\mathcal A\odot \|\theta\|_1$ are convex and $gH$-continuous with $dom\ E=\R^m(\neq \emptyset)$. Therefore, the associated scalarized optimization problem is given by
\begin{equation*}
    \textbf{(SOP5)}\ \min_{\theta\in \R^m}\phi_\lambda(\theta)=\lambda E^l(\theta)+(1-\lambda)E^u(\theta),
\end{equation*}
where $\lambda\in(0,1)$ and $E^l$ and $E^u$ are endpoint functions of $E(\theta)$. By Proposition~\ref{ppsn:2}, the scalarized objective $\phi_\lambda$ is proper, convex and lower semicontinuous, and therefore the proposed PPA-IVOP is directly applicable.\\

\noindent To illustrate the performance of the proposed PPA-IVOP,
we consider the corresponding interval valued dataset reported in Table~\ref{tab:xl}, consisting of $n=12$ interval valued training dataset with $m=2$ interval valued features.

\begin{table}[H]
\centering
\begin{tabular}{cccc}
\toprule
$k$ & $X_k^1$ & $X_k^2$ & $Y_k$ \\
\midrule
1 & $[16.02,\,16.67]$ & $[37.42,\,38.18]$ & $[400.12,\,410.18]$ \\
2 & $[16.58,\,17.01]$ & $[37.96,\,38.52]$ & $[406.34,\,414.22]$  \\
3 & $[16.95,\,17.48]$ & $[38.61,\,39.08]$ & $[414.08,\,421.36]$  \\
4 & $[16.69,\,17.22]$ & $[38.22,\,39.05]$ & $[412.06,\,420.88]$  \\
5 & $[16.47,\,17.02]$ & $[39.01,\,39.66]$ & $[416.12,\,424.91]$  \\
6 & $[16.61,\,16.96]$ & $[39.68,\,40.43]$ & $[422.48,\,431.12]$  \\
7 & $[16.28,\,16.83]$ & $[39.52,\,40.17]$ & $[422.55,\,430.82]$ \\
8 & $[16.91,\,17.18]$ & $[39.89,\,40.31]$ & $[425.43,\,431.17]$   \\
9 & $[16.62,\,17.20]$ & $[40.08,\,40.93]$ & $[426.08,\,437.04]$ \\
10 & $[16.80,\,17.38]$ & $[40.41,\,41.10]$ & $[431.95,\,441.10]$  \\
11 & $[17.05,\,17.71]$ & $[40.47,\,41.18]$ & $[431.88,\,441.02]$ \\
12 & $[16.14,\,16.73]$ & $[41.01,\,41.62]$ & $[434.01,\,440.95]$  \\
\bottomrule
\end{tabular}
\caption{Dataset of the $lasso$ problem.}
\label{tab:xl}
\end{table}
For the numerical implementation, the interval valued regularization parameter is chosen as $\mathcal A=[0.02,0.04]$, the proximal parameter is fixed at $\mu=1$, and the initial iterate is selected as $\theta^0=(8,22)^\top$. The PPA-IVOP is then applied to the (IVOP5) for various choices of $\lambda\in(0,1)$ and the iterative procedure is terminated by the stopping rule $\|\theta^{k+1}-\theta^k\|<10^{-8}.$ The numerical results corresponding to various choices of the scalarization parameter $\lambda$ are summarized in Table~\ref{tab:targets}.

\begin{table}[H]
\centering
\begin{tabular}{cccc}
\toprule
$\lambda$ & Iterations & $\overline{\theta}_1$ & $\overline{ \theta}_2$  \\
\midrule
0.1 & 5 & 0.0000 & 10.6934\\
0.3 & 14 & 0.5531 & 10.4589\\
0.5 & 15 & 1.3683 & 10.1132\\
0.7 & 16 & 2.2079 & 9.7571\\
0.9 & 13 & 3.0733 & 9.3901\\
\bottomrule
\end{tabular}
\caption{Output of PPA-IVOP for (IVOP5).}
\label{tab:targets}
\end{table}
The predicted interval-valued outputs $\mathcal P(X_k,\overline{\theta})$ obtained from the proposed method are compared with the actual target $Y_k$. The comparison table for $\lambda=\{0.1,0.3,0.5,0.7,0.9\}$ are shown in the Table~\ref{tab:estimated_actual}.

\begin{table}[H]
\centering
\scriptsize
\renewcommand{\arraystretch}{1.6}
\setlength{\tabcolsep}{4.5pt}
\begin{tabular}{ccccccc}
\toprule
& \multicolumn{5}{c}{Predicted Output $\mathcal P(X_k,\overline{\theta})$} & \\
\cmidrule(lr){2-6}
$k$ &  $\lambda=0.1$ & $\lambda=0.3$ & $\lambda=0.5$ & $\lambda=0.7$ & $\lambda=0.9$ & Actual $Y_k$ \\
\midrule
1 & [400.15, 408.27] & [400.23, 408.54] & [400.36, 408.93] & [400.48, 409.33] & [400.61, 409.75] & [400.12, 410.18] \\
2 & [405.92, 411.91] & [406.19, 412.29] & [406.58, 412.84] & [406.99, 413.40] & [407.41, 413.99] & [406.34, 414.22] \\
3 & [412.87, 417.90] & [413.19, 418.40] & [413.66, 419.14] & [414.15, 419.90] & [414.65, 420.69] & [414.08, 421.36] \\
4 & [408.70, 417.58] & [408.97, 417.94] & [409.36, 418.48] & [409.77, 419.04] & [410.18, 419.61] & [412.06, 420.88] \\
5 & [417.15, 424.10] & [417.11, 424.25] & [417.05, 424.46] & [416.99, 424.68] & [416.93, 424.91] & [416.12, 424.91] \\
6 & [424.31, 432.33] & [424.20, 432.23] & [424.02, 432.08] & [423.84, 431.93] & [423.65, 431.77] & [422.48, 431.12] \\
7 & [422.60, 429.55] & [422.34, 429.44] & [421.95, 429.28] & [421.55, 429.10] & [421.13, 428.93] & [422.55, 430.82] \\
8 & [426.56, 431.05] & [426.56, 431.10] & [426.55, 431.17] & [426.55, 431.24] & [426.54, 431.32] & [425.43, 431.17] \\
9 & [428.59, 437.68] & [428.39, 437.60] & [428.08, 437.47] & [427.76, 437.34] & [427.44, 437.20] & [426.08, 437.04] \\
10 & [432.12, 439.50] & [431.94, 439.47] & [431.66, 439.43] & [431.38, 439.39] & [431.09, 439.35] & [431.95, 441.10] \\
11 & [432.76, 440.35] & [432.70, 440.49] & [432.61, 440.69] & [432.52, 440.90] & [432.42, 441.11] & [431.88, 441.02] \\
12 & [438.53, 445.06] & [437.85, 444.55] & [436.83, 443.80] & [435.78, 443.03] & [434.69, 442.23] & [434.01, 440.95] \\
\bottomrule
\end{tabular}
\caption{Comparison of the predicted interval outputs.}
\label{tab:estimated_actual}
\end{table}

It can be observed from the above Table~\ref{tab:estimated_actual} that the predicted interval outputs $\mathcal P(X_k;\overline{\theta})$ approximate the actual interval targets $Y_k$ as $\lambda$ sweeps from $(0,1)$, for all $k=1,2,\dots,12$. Furthermore, the solutions corresponding to various choices of $\lambda\in(0,1)$ generates comparable predictions, hence resulting to robustness of the proposed PPA-IVOP on this $lasso$ problem.

\begin{figure}[H]
\centering
\begin{subfigure}{0.42\textwidth}
    \centering
    \includegraphics[width=\linewidth]{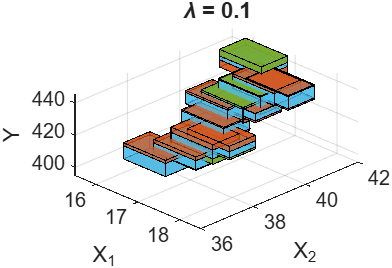}
\end{subfigure}
\hfill
\begin{subfigure}{0.42\textwidth}
    \centering
    \includegraphics[width=\linewidth]{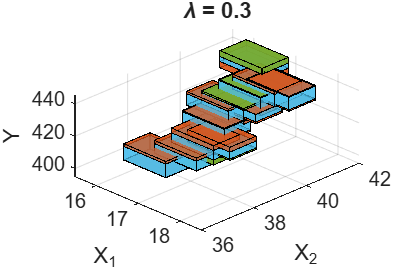}
\end{subfigure}

\vspace{0.5cm}

\begin{subfigure}{0.42\textwidth}
    \centering
    \includegraphics[width=\linewidth]{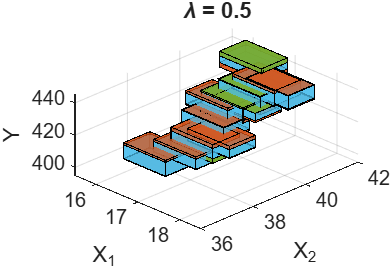}
\end{subfigure}
\hfill
\begin{subfigure}{0.42\textwidth}
    \centering
    \includegraphics[width=\linewidth]{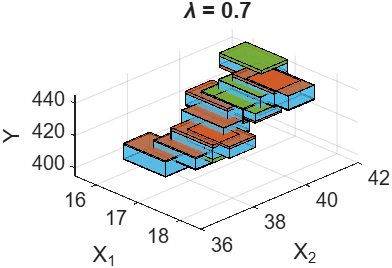}
\end{subfigure}

\vspace{0.5cm}

\begin{subfigure}{0.42\textwidth}
    \centering
    \includegraphics[width=\linewidth]{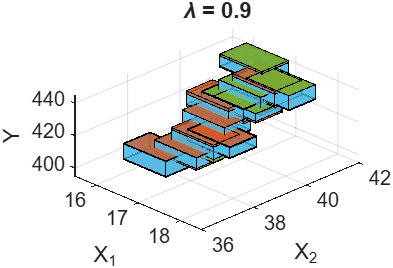}
\end{subfigure}

\caption{Estimated interval outputs with different scalarization parameters.}
\label{fig:lambda_outputs}
\end{figure}

Figure~\ref{fig:lambda_outputs} provides a graphical comparison between the actual outputs and the estimated outputs obtained by the PPA-IVOP for different choices of the scalarization parameter.The light blue, light green and orange regions represents the common region of the estimated and actual target, the extended region of estimated target $\mathcal P(X_k;\bar\theta)$ and the extended region of actual output $Y_k$, respectively.

\section{Conclusion and future research directions}\label{Section Conclusion}
In this paper, we developed a proximal point algorithm applicable to a class of unconstrained interval valued optimization problem (IVOP) whose objective function is a proper convex, and generalized Hukuhara lower semicontinuous. The central idea of the proposed approach is to transform the (IVOP) into an associated scalarized optimization problem (SOP) through a scalarization technique and to establish the relationship between the $LU$-nondominated and optimal solutions between the two problems. Furthermore, we derived the structural properties of the interval valued objective are inherited by the associated scalarized function, thereby guaranteeing the well-posedness of every proximal subproblem. Under the proposed structural assumptions, we established that the PPA-IVOP admits a unique $LU$-nondominated solution of the (IVOP), independent of the choice of the scalarization parameter. Further, the global convergence of the generated sequence to an $LU$-nondominated solution was established, together with a sublinear convergence rate of order $O(\frac{1}{k})$ in the convex case and a linear rate of order $O(\eta^k)$ with $\eta\in(0,1)$ in the strongly convex case. Numerical experiments were presented to demonstrate the results obtained from theoretical convergence analysis of the PPA-IVOP along with a comparative analysis with an existing method, and finally, its applicability was demonstrated through an interval-valued $lasso$ problem, highlighting the potential of the proposed method in the optimization problems involving interval uncertainty.

The proposed framework provides a foundation for several future research directions. Since the scalarization selects LU-nondominated solutions through $\lambda$, alternative scalarization techniques, adaptive choices of $\lambda$, or schemes that work directly in the interval space may allow a more complete exploration of the nondominated set. The analysis can also be extended to inexact and accelerated variants (inertial, variable-metric, quasi-Newton), which would be useful when the proximal operator has no closed form. Finally, one may look at nonconvex, constrained, bilevel and multiobjective interval valued problems, as well as large scale applications with real-world data.

\section*{Author contribution}
All authors contributed equally to this work.

\section*{Funding}
This authors received no external funding for this research.

\section*{Data availability} 
Data will be made available upon reasonable request.

\section*{Conflict of interest}
The authors declare that they have no conflict of interest.

\bibliographystyle{plain}

\end{document}